\documentclass[12pt,a4paper,twoside]{amsart}

\usepackage[T1]{fontenc}
\input xypic\usepackage{amsmath}
\usepackage{amssymb}

\usepackage{amssymb}
\usepackage{latexsym}
\usepackage{graphicx}
\usepackage{psfrag}
\newcommand{\newtext}{}
\newcommand{\newtexx}{}
\newcommand{\newtexxx}{}

\newcommand{\newtexxt}{}

\newcommand{\newwwtext}{}

\newcommand{\R}{{\mathbb R}}

\newcommand{\pM}{{\partial M}}

\newcommand{\bt}{{\bf t}}
\newcommand{\bn}{{\bf n}}
\newcommand{\ud}{{\underline{\delta}}}
\newcommand{\Hnull}{{\stackrel{\circ}H}{}}
\newcommand{\Hbnull}{{\stackrel{\circ}{\bf H}{}}}
\newcommand{\Cnull}{{{\stackrel{\bullet}C}{}^\infty}}

\def\hat{\widehat}
\def\tilde{\widetilde}
\def \bfo {\begin {eqnarray*} }
\def \efo {\end {eqnarray*} }
\def \ba {\begin {eqnarray*} }
\def \ea {\end {eqnarray*} }
\def \beq {\begin {eqnarray}}
\def \eeq {\end {eqnarray}}

\def \det {\hbox{det}}

\def \e {\varepsilon}
\def \p {\partial}
\def \d  {{\delta_{\alpha}}}

\def\M{{\mathcal M}}
\def\F{{\mathcal F}}
\def\cZ{{\mathcal Z}}

\def\hf{{\hat f}}
\def\hh{{\hat h}}
\def\ud{{\underline \delta}}

\newtheorem{definition}{Definition}[section]
\newtheorem{theorem}[definition]{Theorem}
\newtheorem{lemma}[definition]{Lemma}

\newtheorem{proposition}[definition]{Proposition}
\newtheorem{corollary}[definition]{Corollary}
\newtheorem{remark}[definition]{Remark}

\begin{document}
\title
[Inverse Problems for Maxwell's Equations]
{Maxwell's Equations with Scalar Impedance:\\
Inverse Problems with data given on a part of the boundary}

%
%

\author[Yaroslav Kurylev, Matti Lassas and Erkki Somersalo]{
  { Yaroslav Kurylev}\\
  Loughborough University\\
  Department of Mathematical Sciences\\
  Leicestershire LE11 3TU, UK
  \\
  \\
  Matti Lassas\\
  Helsinki University of Technology\\
  Institute of Mathematics\\
  PO Box 1100, 02015 TKK, Finland
  \\
  \\
  Erkki Somersalo\\
  Helsinki University of Technology\\
  Institute of Mathematics\\
  PO Box 1100, 02015 TKK, Finland
}

\makeatletter
\def\@setauthors{%
  \begingroup
  \trivlist
  \centering \@topsep30\p@\relax
  \advance\@topsep by -\baselineskip
  \item\relax
  \andify\authors
  \def\\{\protect\linebreak}%
{\authors}%
  \endtrivlist
  \endgroup
}
\makeatother

\maketitle


{\bf Abstract:} We study Maxwell's equations in time domain in an
anisotropic medium. The goal of the paper is to solve an inverse
boundary value problem for anisotropies characterized by scalar impedance $\alpha$.
This means that the material is conformal, i.e., the electric permittivity
$\epsilon$ and magnetic permeability $\mu$ are tensors satisfying
$\mu =\alpha^2\epsilon$. This condition is equivalent
to a single propagation speed of waves with different polarizations
which uniquely defines an underlying Riemannian structure.
The analysis is based on an invariant formulation
of the system of electrodynamics as a Dirac type first order system
on a Riemannian $3-$manifold with an additional structure of the wave
impedance, $(M,g,\alpha)$, where $g$ is the travel-time metric.
We study the properties of this system in the first part of the paper.
In the second part we consider the inverse  problem, that is,
the determination of
$(M,g,\alpha)$ from measurements done only on an open part 
of the boundary and on a finite time interval.
As an application, in the isotropic case with $M\subset \R^3$, 
we prove that the boundary data given only on an open part of the boundary
determine uniquely  the domain $M$ and the coefficients $\epsilon$ 
and $\mu$.


\noindent
{\bf Keywords:} Maxwell's equations,  manifolds, inverse problems. 


\section*{Introduction}\label{introduction}

In this paper we study direct and inverse 
boundary value problems for Maxwell's equations
on Euclidean domains in $\R^3$ and on  compact  manifolds.
In a bounded smooth domain $M^{\rm int} \subset \R^3$,
Maxwell's equations  for  the electric
and magnetic fields $E$ and $H$ and the associated electric
displacement $D$ and magnetic flux density $B$ are
\beq
\label{Maxwell-Faraday 1}
& & {\rm curl}\,E(x,t) = - B_t(x,t),\mbox{ (Maxwell--Faraday)},\quad
{\rm div}\,B(x,t)=0,\\
\label{Maxwell--Ampere 1}
& & {\rm curl }\,H(x,t) =\phantom{-} D_t(x,t),
\mbox{ (Maxwell--Amp\`{e}re)}, \quad
   {\rm div}\,D(x,t)=0.
\eeq
Under the assumption of non-dispersive, linear and non-chiral medium,
these are augmented with the constitutive relations
\beq
\label{constitutive 1}
  D(x,t) = \epsilon(x)E(x,t),\quad
  B(x,t) = \mu(x)H(x,t).
\eeq
Here electric permittivity $\epsilon$ and magnetic permeability $\mu$
are smooth $3 \times 3$  time-independent positive matrices.

The present work consists of two parts.
In the first part, we pursue further the invariant formulation of
Maxwell's equations (\ref{Maxwell-Faraday 1})--(\ref{constitutive 1})
for certain anisotropic materials,
characterized by a {\em scalar wave impedance}. This means that
the material is conformal in the sense that
\beq
\label{impedance 1}
  \mu =\alpha^2\epsilon,
\eeq
where $\alpha$ is a positive scalar function.
{\newwwtext   Physically, the condition of a scalar wave impedance is tantamount to
a single propagation speed of waves with different polarizations
and,
 without it, the wave propagation
gives rise to several  generally non-Riemannian metrics, see  e.g. \cite{Ka2}.
Therefore, we consider direct (and inverse) problems for the most 
general subclass of Maxwell's equations which is distinguished by the fact 
that electromagnetic fields with different polarization propagate with the same 
velocity which, of cause, may
depend on the propagation direction.}
The scalar wave impedance
is  encountered in many physical situations.
For instance,
in a curved spacetime with coordinates
$(x,t)\in \R^3\times \R$ and a ``time-independent''
metric $ds^2=g_{jk}(x)dx^jdx^k-dt^2$,
Maxwell's equations with scalar permittivity and permeability
correspond in the coordinate invariant form to Maxwell's equations
(\ref{Maxwell-Faraday 1})--(\ref{constitutive 1})
with scalar wave impedance,
see \cite [Sec. 14.1.c]{Frankel} or \cite[Sec. 90]{LanL}.
{\newtexxx Clearly, all isotropic media, i.e. with scalar 
$\epsilon$ and $\mu$ have scalar wave impedance.}

The invariant approach leads us to formulate Maxwell's equations on
3--manifolds as a first order Dirac type system. From the operator
theoretic point of view, this formulation is based on an
elliptization procedure by extending Maxwell's equations to the
bundle of exterior differential forms over the manifold. This is a generalization of
the elliptization of Birman
and Solomyak
and Picard (see \cite{Birman,picard}).

In the second part of the work, we consider
the inverse boundary value problem for
Maxwell's equations with scalar wave impedance.
In physical terms,
the goal is to determine the material parameter tensors $\epsilon$ and
$\mu$ in a bounded domain from field observations at the boundary 
{\newtexxx or a part of the boundary} of that
domain.

In the invariant approach to  Maxwell's
equations, the domain $M$ is considered as a 3--manifold
and the vector fields $E$, $H$, $D$ and $B$ as
differential forms. This
alternative formulation has several advantages both
from the theoretical and practical points of view.
First, the invariance of the system and the boundary
measurements 
{\newtexxx with respect to 
diffeomorphisms of $M$ that preserve the part of the boundary
where these measurement are done} 
is essential for the inverse problem.
It is possible to prove unique identifiability in the invariant formulation
and then use this result to completely
characterize group of transformations between indistinguishable
parameters $\epsilon$ and $\mu$ in the case $M\subset \R^3$.
{\newtexxx In particular, when $\epsilon$ and $\mu$ 
are scalar functions, this result implies uniqueness of the
determination of $\epsilon$, $\mu$, and $M$ from data 
on a part of boundary.}  

Second, the formulation of electrodynamics in terms of
differential forms reflects the way in which these
fields are actually observed. For instance, flux
quantities are expressed as 2--forms while field
quantities that correspond to forces are naturally
written as 1--forms.
This point of view has been adopted in modern physics,
at least in the free space, see
e.g. \cite{Frankel}, as well as in applications where the numerical
treatment of the equations is done using the Whitney elements.
An extensive treatment of this topic can be found in
\cite{bossavit1,bossavit2}.
For the original reference concerning the Whitney elements,
see \cite{whitney}.

As inverse problems of electrodynamics have a great significance
in physics and applications,  they have been studied starting from the 30's,
see e.g. \cite{langer,schlichter}, where the one-dimensional
case was considered. However, results concerning
the multidimensional inverse problems in electrodynamics are relatively recent.
The first breakthrough achieved in   \cite{SIC,CP,OPS,OS} was based on the
use of complex geometrical optics.  
{\newwwtext These papers were devoted to the identifiability of isotropic 
material parameters $\epsilon$ and $\mu$ from the fixed-frequency  data collected on 
$\p M$, namely, the stationary admittance map (for the definition of the time-dependent
admittance map see Definition \ref{27.11.d}  in section 1.3.) Under some mild 
geometric assumptions it was shown there that these data determine isotropic 
$\epsilon$ and $\mu$ and also isotropic conductivity, $\sigma$, uniquely.}
These works were based on the ideas previously developed in
\cite{SyU,Na1,Na2} to tackle the
scalar Calder\'{o}n problem, introduced in \cite{Cl}.
Other approaches to the isotropic inverse problem  for Maxwell's
equations work directly in the time domain, see \cite{Be3,rom}.
{\newwwtext Regarding the case $\sigma=0$ which is considered in this paper,
the result obtained in \cite{rom} proves the identifiability of  $\epsilon$ and $\mu$
from the time-dependent data collected on the whole $\p M$ in the case when
$M$, considered as a Riemannian manifold with metric 
$dl^2=\epsilon \mu |dx|^2$, is  simple geodesic. Constructions of \cite{Be3} make it possible
to find the product, $\epsilon \mu$ of unknown parameters $\epsilon, \, \mu$.
Moreover, the results of \cite{Be3} are of a local nature making it possible to find this product
only in some collar neighbourhood of $\p M$. Time-dependent inverse problem for
isotropic Maxwell's equations was also considered in \cite{BI} which used a time 
Fourier transform to reduce the problem to the one in the frequency domain so
that to apply the results of \cite{OPS,OS}.}



Much less is known in the anisotropic case, where material
parameters are matrix valued functions.
{\newwwtext  The case of anisotropic $\epsilon=\mu$ was considered in \cite{BeIsPSh} 
where it was
shown that the time-dependent admittance map known on $\p M$ makes it possible to recover $\epsilon=\mu$ locally, i.e. in some collar neighbourhood of $\p M$.
  In spite of very little knowledge, it is, however,}
  clear from the study of
scalar anisotropic problems  that, instead of uniqueness, one obtains
uniqueness only up to a group of transformations, involving proper
coordinate changes, see e.g. \cite{LeU,Sy,
BeKu3,Ku1, KK2, LU,LTU}.
A similar result for Maxwell's equations was conjectured in \cite{Sy2},
based on the analysis of the linearized inverse problem.
Therefore, it is natural to split the study of this problem into
two steps.
First, to formulate and solve the corresponding coordinate-invariant
inverse problem, i.e., an  inverse problem on a manifold. Second,
to analyse the properties resulting from an embedding of the
manifold into $\R^3$. For a systematic development of this
approach, see \cite{KKL}.

In recent years,  inverse problems with data on a part of the 
boundary have attracted much interest, see \cite{GU,KK2,KSU,IU,LU}.
 Part of the motivation
come from the physical setting when only a part of the boundary
is accessible. 
{\newwwtext However, as far as we know, there are currently no results on 
identifiability of the shape of the domain $M$ and/or the material 
parameters $\epsilon,\,\mu$ on it from inverse data collected on an arbitrary open
subset, $\Gamma \subset \p M$.}

A fruitful approach to scalar inverse problems,
including those with data on a part of the boundary,
turned out to be the boundary control method,
originated in \cite{Be1} for the isotropic acoustic wave equation.
In the  anisotropic context,  it has been developed
for the Laplacian on Riemannian manifolds \cite{BeKu3} and for general
anisotropic self-adjoint \cite{Ku2,Ku3} ) and certain
non-selfadjoint inverse problems \cite{KL1}.


The current article pursues the study 
{\newwwtext of inverse problems for Maxwell's equations  significantly further dealing with the global reconstruction of the shape of the domain  or, more general, 
$3-$manifold $M$,   metric tensor $g$ and scalar wave impedance $\alpha$, the 
latter two being equivalent to the reconstruction of  $\epsilon$ and $\mu$. Being based
on the boundary control, the method developed here 
 combines ideas of the articles
\cite{KL2} and \cite{KL3}
with those of \cite{OPS} and \cite{OS}. What is more, to be able to study
anisotropic Maxwell's equations, we introduce two essential new ideas.
First, we characterize the subspaces controlled from the boundary {\em by
duality}, thus avoiding the difficulties arising from the complicated
topology
of the domains 
of influence but still providing necessary information about the structure of
achievable sets, see e.g. Theorem \ref{local control th} in section 1.5. This makes our approach 
much different from that in \cite{BeIsPSh,Be3}  also based
on the boundary control method. Indeed, the method of \cite{BeIsPSh,Be3} 
requires local controllability
in the domains of influence which is no more valid for large times,
e.g. \cite{Be-Gl2} thus making constructions of \cite{BeIsPSh,Be3} inappropriate outside a collar
neighbourhood  of $\p M$. Second,  we develop a {\em method of focused waves} which
enables us to recover pointwise values of electromagnetic
waves
 on the manifold and, therefore, reconstruct not
only the metric $g$, as in \cite{Be3}, but also the impedance $\alpha$.
}


The main results of this paper can be summarized as follows.
\begin{enumerate}

\item The knowledge of the complete dynamical boundary data 
over a  sufficiently large
finite period of time determines uniquely
the compact manifold endowed with the travel time metric
as well as the scalar wave impedance (Theorem \ref{ip}).
This is valid also when measurements are made on
a part of the boundary (Theorem \ref{ip2}).
{\newtexxx The necessary time of observation is double of the time
required to fill the manifold from the observed part of the boundary. 
}

\item For the corresponding anisotropic
inverse boundary value problem with scalar wave impedance for
bounded domains in $\R^3$, the non-uniqueness 
is completely characterized by describing
the class of possible transformations between material tensors
that are indistinguishable from the observed part of the boundary (Theorem \ref{group}).

\item
{\newtexx For the corresponding isotropic
inverse boundary value problem for
bounded domains in $\R^3$, the shape of the domain and  
the material
parameters inside it} are uniquely determined from measurements
done on a part of the boundary (Theorem \ref{uniqueness in isotropic case}).

\end{enumerate} 
{\newtexx
Some of the results of the paper have been announced in \cite{KLS,KLS2}.
}

\section{Maxwell's equations on a manifold}

This chapter is devoted to Maxwell's equations on a compact
oriented  3--manifold with boundary. We concentrate on the properties
of these equations important for the inverse problem considered
in Chapter 2.

{\newtext We start with the formulation of Maxwell's equations for 1--
and 2--forms. These equations are augmented to the complete Maxwell
  system on the full
bundle of exterior differential forms over a 3--dimensional Riemannian manifold.}
This allows us to define and analyze properties
of an elliptic operator related to  Maxwell's equations and to
study the corresponding initial boundary value problem.
Crucial results of Sections 1.3 and 1.4 are the Blagovestchenskii
formula, Theorem \ref{blacho}, enabling us to
evaluate inner products of electromagnetic waves in terms
of the admittance map ${\mathcal Z}$, also defined in Section 1.3, and
the unique continuation result for  Maxwell's equations with Cauchy data
on the lateral boundary. Building on these results, we obtain
local and global controllability for electromagnetic waves
generated by boundary sources and define, in a usual manner, spaces
of generalized sources.

\subsection{Invariant definition of Maxwell's equations}

To define Maxwell's equations invariantly, consider
a smooth compact oriented connected Riemannian 3-manifold $M$, $\partial M\neq
\emptyset$,
with a metric $g_0$, that we call the background metric.
Clearly, in physical applications we take $M\subset \R^3$ with $g_0$
being the Euclidean metric.
Analogously to (\ref {Maxwell-Faraday 1}) and (\ref{Maxwell--Ampere 1}),
Maxwell's equations on the manifold $M$
are  equations of the form
\beq
\label{Maxwell-Faraday}
& & {\rm curl}\,E(x,t) = - B_t(x,t),
\quad {\rm div}\,B(x,t)=0,\\
\label{Maxwell--Ampere}
& & {\rm curl }\,H(x,t) =\phantom{-} D_t(x,t), \quad
   {\rm div}\,D(x,t)=0.
\eeq
Here  $E$, $H$, $D$, $B \in \Gamma M$, the space
of $C^\infty$-smooth vector fields on $M$. They are
related by the constitutive relations,
\beq
\label{constitutive}
  D(x,t) = \epsilon(x)E(x,t),\quad
  B(x,t) = \mu(x)H(x,t),
\eeq
where $\epsilon$ and $\mu$ are $C^\infty$-smooth positive
definite $(1,1)$--tensor fields on $M$. We remind that, for $X\in\Gamma M$,
\beq \label{A 23}
({\rm curl}\,X)^\flat = *_0 dX^\flat,\quad
{\rm div}\,X = -*_0 d *_0 X^\flat.
\eeq
Here, $d$ is the exterior differential,
{\newwwtext $^\flat$ is the fiberwise duality  between 1--forms and vector fields
\bfo
X \in \Gamma M \rightarrow X^\flat \in \Omega^1M, \quad
X^\flat(Y)= g_0(X,Y),
\efo
with $\Omega^1M$ and, generally, $\Omega^kM$ standing for the bundle of
differential $k-$forms on $M$. At last $*_0$ is the Hodge operator with respect
to metric $g_0$, acting fiberwise,
\bfo
*_0:\Omega^kM \rightarrow \Omega^{3-k}M.
\efo
Throughout the paper, we assume that the
wave impedance is scalar, i.e.,
\beq
\label{impedance}
  \mu =\alpha^2\epsilon,
\eeq
where $\alpha=\alpha(x)$ is a positive  scalar function. 
This allows us to introduce a new metric, $g$ on $M$.
Considering it as a quadratic form  on differential 1-forms,
we have
\beq
\label{metric}
g(X^\flat,Y^\flat) = \frac 1{{\rm det}\, z}g_0(X^\flat,zY^\flat), \quad
g^{ij}= \frac{1}{\det{z}}g_0^{ik}z^j_k,
\quad z=\alpha\epsilon =\frac
1\alpha\mu,
\eeq
As we see later, this metric is responsible for the velocity of electromagnetic 
wave propagation and we call it the {\it travel-time metric}.
}


{\newwwtext  Introduce the differential  $1-$ and $2-$forms, $\omega^1, \,\nu^1$ and 
$\omega^2, \,\nu^2$,
\beq
\label{forms}
\omega^1= E^\flat,\,\, \nu^1= \alpha H^\flat,
\quad \omega^2= *_0B^\flat,\,\, \nu^2= *_0\alpha D^\flat.
\eeq
Applying the operator $^\flat$ to Maxwell's equations (5)-(6), we rewrite them 
in terms of $\omega^1, \,\nu^1,\, \omega^2, \,\nu^2$,	
\beq
\label{additional}
& &
d \omega^1= -\omega^2_t, \quad d \omega^2=0,
\\
\nonumber
& &
d \frac{1}{\alpha} \nu^1 = \frac{1}{\alpha} \nu^2_t, \quad d \frac{1}{\alpha} \nu^2=0.
\eeq
Now, using (\ref{metric}) it is a straightforward matter to check that the constitutive relations 
(\ref{constitutive})
assume the form,
\beq
\label{consitutive1}
\nu^2= * \omega^1, \quad
\omega^2 = * \nu^1,
\eeq
where $*$ is the Hodge operator with respect to the travel-time metric, $g$.
Eliminating, by means of (\ref{consitutive1}), $\nu^{1,2}$ from equations
(\ref{additional}), we transform (\ref{Maxwell-Faraday})-(\ref{constitutive}) into the system
\beq
\label{MF}
& &\omega^1_t = \delta_{\alpha} \omega^2, \quad \delta_{\alpha} \omega^1=0,\\
& &
\label{MA}
\omega^2_t = -d \omega^1, \quad d \omega^2=0.
\eeq
Here,  $\delta_{\alpha} :  \Omega^{k}M\to \Omega^{3-k}M$ is the
$\alpha$--codifferential, given by
\beq
\label{codifferential}
\delta_{\alpha} \omega^k = (-1)^k *\alpha \,d \frac{1}{\alpha} \, * \,\omega^k,
\eeq
and $*$ is the Hodge operator in the metric $g$.
These equations are called {\it Maxwell's equations for forms}
on the Riemannian manifold with a scalar
wave impedance $(M,g,\alpha)$.
}

To extend the above equations  to the full bundle of exterior differential forms
${\bf \Omega}M =\Omega^0 M\times\Omega^1 M\times\Omega^2 M
\times\Omega^3 M$, we introduce
  auxiliary  forms, $\omega^0\in\Omega^0 M$ and $\omega^3\in\Omega^3 M$,
which vanish in the electromagnetic theory, by
\[
  \omega_t^0 = \delta_\alpha\omega^1,\quad
  \omega^3_t=-d\omega^2.
\]
Since $\omega^0=0$ and $\omega^3=0$ in electromagnetics,
we can modify  equations (\ref{MF}) and (\ref{MA})  to read
\beq\label{MA1}
& &\omega_t^1=- d\omega_0 +\delta_\alpha\omega^2,
\quad \omega^3_t=-d \omega^2,\\
\label{MF1}
& &\omega_t^2=- d\omega^1 +\delta_\alpha\omega^3,
\quad \omega^0_t=\delta_{\alpha} \omega^1,
\eeq
or, in the matrix
form,
\begin{equation}\label{complete}
  \omega_t + {\mathcal M}\omega=0,
\end{equation}
where
$
  \omega = (\omega^0,\omega^1,\omega^2,\omega^3) \in {\bf \Omega}M,
$
and the operator ${\mathcal M}$ (without prescribing its domain at this
point, i.e., defined as a differential expression)
is given as
\begin{equation}\label{M}
{\mathcal M} = \left(\begin{array}{cccc}
0 &-\delta_{\alpha} &0       &0       \\
d & 0      &-\delta_{\alpha} &0       \\
0 & d      &0       &-\delta_{\alpha} \\
0 & 0      &d       &0
\end{array}\right).
\end{equation}
Equations (\ref{complete}),  (\ref{M}) are called {\em
the complete Maxwell system}.

Note that, identifying
${\bf \Omega}M$ with $\Omega^0M \oplus  \Omega^1M \oplus
  \Omega^2M \oplus \Omega^3M$,  the
complete Maxwell operator can be written as
\beq
\label{Dirac}
{\mathcal M} = d -\d,
\eeq
i.e., it becomes a Dirac type operator on ${\bf \Omega}M$.

An important property of ${\mathcal M}$ is that
\[
  {\mathcal M}^2 = -{\rm diag}(\Delta_\alpha^0,
\Delta_\alpha^1,\Delta_\alpha^2,\Delta_\alpha^3)
  =-{\bf \Delta}_\alpha,
\]
where the operator $\Delta_\alpha^k$ acts on  $k$--forms as
\beq
\label{Laplace}
  \Delta_\alpha^k = d\delta_\alpha + \delta_\alpha d
  =\Delta^{k}_g + Q^k(x,D).
\eeq
Here $\Delta_g^k$ is the Laplace-Beltrami operator
in the metric $g$ and $Q^k(x,D)$
is a first order perturbation.
Hence, if $\omega$ satisfies  equation
(\ref{complete}), it satisfies also wave equation
\begin{equation}\label{wave equation}
  (\partial_t^2 + {\bf \Delta}_\alpha)\omega
  =(\partial_t-{\mathcal M})(\partial_t
+{\mathcal M})\omega =0.
\end{equation}
This formula legitimates the notion of the travel time metric
and makes clear that in the Maxwell system with scalar impedance,
the electromagnetic waves of different
polarization propagate with the same speed determined by the metric $g$.

We end this section with representation of the energy of  electric  and
magnetic fields in terms of
the corresponding differential forms, setting
\begin{eqnarray*}
  {\mathcal E}(E) &=& \frac 12 \int_M g_0(\epsilon E, E)\, dV_0
=\frac 12 \int_M\frac 1\alpha \omega^1\wedge
  *\omega^1,\\
\noalign{\vskip4pt}
  {\mathcal E}(B) &=& \frac 12 \int_M g_0(\mu H, H)\, dV_0
= \frac 12 \int_M\frac 1\alpha \omega^2\wedge
  *\omega^2,
\end{eqnarray*}
where $dV_0$ is volume form of $(M,g_0)$.
These formulae serve as a motivation for our definition
of  inner products in the following section.

\subsection{The Maxwell operator}

In this section we establish a number of notational
conventions and definitions concerning the differential
forms used in this paper.

We define the $L^2$--inner products for $k$--forms
in $\Omega^k M$ as
\[
  (\omega^k,\eta^k)_{L^2} = \int_M \frac 1\alpha
  \omega^k\wedge * \eta^k,\quad
  \omega^k,\;\eta^k\in\Omega^k M,
\]
and denote by
  $L^2(\Omega^k M)$ the completion of
$\Omega^k M$ in the corresponding norm.
We also define
\[
  {\bf L}^2(M) =
  L^2(\Omega^0 M)\times
L^2(\Omega^1 M)\times
L^2(\Omega^2 M)\times L^2(\Omega^3 M),
\]
with Sobolev spaces ${\bf H}^s(M)$,  ${\bf H}_0^s(M)$, $s \in \Bbb{R}$,
given as
\[
{\bf H}^s(M) = H^s(\Omega^0M)\times
H^s(\Omega^1 M)\times
H^s(\Omega^2 M)\times H^s(\Omega^3 M),
\]
\[
{\bf H}_0^s(M) = H^s_0(\Omega^0M)\times
H^s_0(\Omega^1 M)\times
H^s_0(\Omega^2 M)\times H^s_0(\Omega^3 M).
\]
Here, $H^s(\Omega^kM)$ is the Sobolev space of $k$--forms
and $H^s_0(\Omega^k M)$ is the closure in $H^s(\Omega^k M)$
of  the set of $k$--forms in $\Omega^k M$,
which vanish near $\p M$.

Clearly, the natural domain of the exterior derivative, $d$ in
$L^2(\Omega^k M)$ is
\[
  H(d,\Omega^k M) = \left\{ \omega^k\in L^2
  (\Omega^k M)\mid
  d\omega^k \in L^2(\Omega^{k+1} M)\right\},
\]
and the natural domain of  $\d$ is
\[
  H(\delta_\alpha,\Omega^k M) = \left\{ \omega^k
  \in L^2(\Omega^k M)\mid
  \delta_\alpha\omega^k \in L^2(\Omega^{k-1} M)\right\}.
\]
In the sequel, we  drop the sub-index $\alpha$
from the codifferential.

The operators $d$ and $\delta$ are adjoint, i.e. for $C^\infty
_0$--forms $\omega^k, \eta^{k+1}$,
\[
(d\omega^k,\eta^{k+1})_{L^2} = (\omega^k,\delta\eta^{k+1})_{L^2}.
\]
To extend this formula to less regular forms,
let us fix some notations. For $\omega^k
\in\Omega^k M$, we define its {\em tangential}
and {\em normal} boundary components on $\partial M$ as
\begin{eqnarray*}
\bt\omega^k = i^*\omega^k,\quad
\bn\omega^k = i^*(\alpha^{-1}*\omega^k),
\end{eqnarray*}
respectively, where $i^*:\Omega^k M\to \Omega^k\partial M$ is the
pull-back of the natural imbedding  $i:\partial M\to M$. 
With these notations, Stokes' formula for forms can be written as
\begin{equation}\label{stokes1}
  (d\omega^k,\eta^{k+1})_{L^2}-
(\omega^k,\delta
  \eta^{k+1})_{L^2} =
  \langle\bt \omega^k,
\bn \eta^{k+1}\rangle,
\end{equation}
where,
for $\omega^k\in\Omega^kM$ and $\eta^{k+1}\in
\Omega^{k+1}M$,
\[
\langle\bt\omega^k,
\bn\eta^{k+1}\rangle =
   \int_{\partial M}\bt \omega^k\wedge\bn \eta^{k+1}.
\]
There are well defined  extensions of the boundary
trace operators $\bt$ and $\bn$ to $H(d,\Omega^k M)$
and $H(\delta,\Omega^k M)$.
The following result is due to Paquet \cite{paquet}:

\begin{proposition}\label{paquet}
The operators $\bt$ and $\bn$
can be extended to continuous surjective maps
\begin{eqnarray*}
\bt: H(d,\Omega^k M) \to
H^{-1/2}(d,\Omega^k\partial M), \quad\quad
\bn: H(\delta,\Omega^{k+1} M) \to
H^{-1/2}(d,\Omega^{2-k}\partial M),
\end{eqnarray*}
where the space $H^{-1/2}(d,\Omega^k\partial M)$ is the
space of $k$-forms $\omega^k$
on $\partial M$ satisfying
\[
  \omega^k\in H^{-1/2}(\Omega^k \partial M),\quad
  d\omega^k\in H^{-1/2}(\Omega^{k+1} \partial M).
\]
\end{proposition}

Formula (\ref{stokes1}) is instrumental for characterizing the
spaces of forms with vanishing
boundary data. Introducing  $\Hnull(d,\Omega^k M)=\hbox{Ker}\,(\bt)$
and
$\Hnull(\delta,\Omega^{k+1}M)=\hbox{Ker}\,
(\bn)$ and applying Stokes' formula, one can prove in standard way
the following lemma.

\begin{lemma}\label{adjoints}
The adjoint of the operator
\[
  d:L^2(\Omega^k M)\supset H(d,\Omega^k M)\to L^2(\Omega^{k+1} M)
\]
is the operator
$
  \delta:L^2(\Omega^{k+1}M)
\supset\Hnull(\delta,\Omega^{k+1}M)
\to L^2(\Omega^kM)
$
and {\em vice versa}. Similarly, the adjoint of
\[
  \delta:L^2(\Omega^{k+1}M)
\supset H(\delta,\Omega^{k+1}M)
\to L^2(\Omega^kM)
\]
is the operator
$
  d:L^2(\Omega^k M)\supset \Hnull(d,\Omega^k M)\to L^2(\Omega^{k+1} M).
$
\end{lemma}

When there is no risk of confusion we will write for brevity
$H(d)=H(d,\Omega^k M)$ and similarly, {\em mutatis mutandis}
for the other spaces.

For later reference, we point out that
Stokes' formula for the complete Maxwell system
can be written  as
\begin{equation}\label{stokes for system}
  (\eta,{\mathcal M}\omega)_{{\bf L}^2} +  ({\mathcal M} \eta,
\omega)_{{\bf L}^2} =
\langle\bt\omega,\bn\eta\rangle +
\langle\bt\eta,\bn\omega\rangle,
\end{equation}
where $\omega \in{\bf H}$ with
\beq
\label{H}
{\bf H} = H(d)\times[H(d)\cap
H(\delta)]\times[H(d)\cap
H(\delta)]\times H(\delta),
\eeq
$\eta \in {\bf H}^1(M)$,
$
\bt\omega = (\bt\omega^0,\bt\omega^1,\bt\omega^2)$,
$
\bn\omega = (\bn \omega^3,\bn \omega^2,\bn
  \omega^1),
$
  and
\[
\langle\bt\omega,\bn\eta\rangle =
\langle\bt\omega^0,
\bn\eta^1\rangle + \langle\bt\omega^1,
\bn\eta^2\rangle +\langle\bt\omega^2,
\bn\eta^3\rangle.
\]

With these notations, we give the following definition
of the Maxwell operators with electric
boundary condition.

\begin{definition}\label{d. 2}
The Maxwell operator with
electric boundary condition, ${\mathcal M}_{\rm e}$, is an operator in
${\bf L}^2(M)$, with
\[
  {\mathcal D}({\mathcal M}_{\rm e})=\Hbnull_{\hspace{-1mm}\bt}:= \Hnull(d)
\times[\Hnull(d)\cap H(\delta)]\times
  [\Hnull(d)\cap H(\delta)]
\times H(\delta),
\]
and ${\mathcal M}_{\rm e}\omega$, $\omega \in {\mathcal D}({\mathcal 
M}_{\rm e})$ is given
by the differential expression
(\ref{M}).
\end{definition}

In terms of physics, the electric boundary condition is associated with
electrically
perfectly conducting boundaries, i.e., $n\times E =0$, $n\cdot B=0$, where $n$
is the exterior normal vector at the boundary. In
terms of differential forms, this means simply that
$\bt E^\flat=\bt\omega^1 = 0$ and $\bt *_0 B^\flat=
\bt\omega^2 = 0$.
Although not used in the sequel, the Maxwell operator with
magnetic boundary condition, ${\mathcal M}_{\rm m}$, is given by (\ref{M}) with
the domain
\[
  {\mathcal D}({\mathcal M}_{\rm m})=\Hbnull_{\hspace{-1mm}\bn}:= H(d)
\times[H(d)\cap \Hnull(\delta)]\times
  [H(d)\cap \Hnull(\delta)]
\times \Hnull(\delta).
\]

Consider the intersections of spaces
in the  definition of ${\mathcal D}({\mathcal M}_{\rm e})$ and
${\mathcal D}({\mathcal M}_{\rm m})$.
Let
\bfo
& & \Hnull^1_\bt(\Omega^k M)=\{\omega^k\in H^1
(\Omega^k M)\mid \bt\omega^k=0\}, \\
& &\Hnull^1_\bn(\Omega^k M)=\{\omega^k\in H^1
(\Omega^k M)\mid \bn\omega^k=0\}.
\efo
It is a direct consequence of Gaffney's inequality
(see \cite{Sc}) that
\beq
\nonumber& & \Hnull(d,\Omega^k M)\cap H(\delta,\Omega^k M)
  =\Hnull^1_\bt(\Omega^k M), \\
\nonumber
& & H(d,\Omega^k M)\cap \Hnull(\delta,\Omega^k M)
  =\Hnull^1_\bn(\Omega^k M).
\eeq

The following lemma is a straightforward application
of Lemma
\ref{adjoints} and the classical
Hodge-Weyl decomposition \cite{Sc}.

\begin{lemma}\label{lem 2.4}
  The electric Maxwell operator has the following properties:
\begin{enumerate}

\item[(i)] The operator ${\mathcal M}_{\rm e}$ is skew-adjoint.
\item[(ii)] The operator ${\mathcal M}_{\rm e}$ defines an elliptic
differential operator in $M^{\rm int}.$
\item[(iii)] ${\rm Ker}\,({\mathcal M}_{\rm e})=\big\{(0,\omega^1,\omega^2,
\omega^3)\in \Hbnull_{\hspace{-1mm}\bf t}\mid\
d\omega^1=0,\ \delta\omega^1=0,\
d\omega^2=0,\ \delta\omega^2=0,\ \delta\omega^3=0\big\}$.

\item[(iv)] ${\rm Ran}\,({\mathcal M}_{\rm e}) =
  L^2(\Omega ^0M)\times \big(\delta H(\delta,\Omega^2 M)
+d\Hnull(d,\Omega^0 M)\big)\times$ \\
$\hbox{ } \quad \quad \quad \quad \quad \times \big(\delta H(\delta,\Omega^3 M)
+d\Hnull (d,\Omega^1 M)\big)\times d\Hnull (d,\Omega^2 M)$.
\end{enumerate}
\end{lemma}

By the skew-adjointness, it is possible to define
weak solutions to initial boundary-value problems
needed later.

\subsection{Initial--boundary value problem}

In the sequel, we denote the
forms $\omega(x,t)$ by $\omega(t)$ or $\omega$
when there is no danger of misunderstanding.

By the {\em weak solution} to the initial boundary value
problem
\beq\label{ibvp}
& &\p_t \omega + {\mathcal M}\omega =\rho
\in  L^1_{\rm loc}(\R,{\bf L}^2(M)),
\\
& &\nonumber
\bt\omega\big|_{\partial M\times \R}=0,\quad
\omega(0)=\omega_0\in{\bf L}^2,
\eeq
we mean the form $\omega(t) \in C(\R,{\bf L}^2(M))$ defined as
\beq\label{weak solution}
\omega(t) ={\mathcal U}(t)\omega_0+\int_0^t{\mathcal U}
(t-s)\rho(s)ds,
\eeq
where ${\mathcal U}(t)={\rm exp}(-t{\mathcal M}_
{\rm e})$ is the unitary operator in ${\bf L}^2$ generated by
${\mathcal M}_{\rm e}$.
Similarly, we define weak solutions with initial
data at $t=T$,  $T\in \R$.
Assuming $\rho \in C(\R,{\bf L}^2)$, the solution has
more regularity, $\omega\in C(\R,{\bf L}^2)\cap C^1(\R,{\bf H}')$,
where ${\bf H}'$ denotes the dual of ${\bf H}$.

We consider also the normal boundary trace
$\bn\omega$ of
weak solutions. Omitting the details, we note that
$\bn\omega$ is defined as a limit of smooth solutions
approximating the weak one.

The following result gives a sufficient condition for a weak solution of
the complete system to be also a solution of Maxwell's equations.

\begin{lemma}\label{weak is maxwell}
Assume that the initial data $\omega_0$ is of the
form $\omega_0=(0,\omega_0^1,\omega_0^2,0)$, where
$\delta\omega_0^1 =0$, $d\omega_0^2=0$, and
$\rho=0$. Then the weak solution $\omega(t)$ of the form
(\ref{weak solution}) satisfies also
Maxwell's equations (\ref{MF}),  (\ref{MA}),
i.e., $\omega^0=0$ and  $\omega^3=0$.
\end{lemma}

{\em Proof:}
As seen from (\ref{wave equation}),
  $\omega^0(t)$ satisfies
the wave equation
\[
  \Delta^0_{\alpha}\omega^0 + \omega^0_{tt}=0,
\]
  with the Dirichlet
boundary condition $\bt\omega^0 =0$. The initial
data for $\omega^0$ is
\ba
& & \omega^0(0) = \omega_0^0 =0,\quad
  \omega^0_t(0) = \delta\omega^1\big|_{t=0}
  = \delta\omega^1_0=0.
\ea
Hence,  $\omega^0(t)=0$ for all $t \in \R$.

Similarly, $\omega^3(t)$ satisfies the wave equation
with the initial data
\ba
\omega^3(0) = \omega_3^0 =0,\quad
  \omega^3_t(0) = -d\omega^2\big|_{t=0}
  = -d\omega^2_0=0.
\ea
As for the boundary condition, we observe that
\[
  \bt \delta\omega^3
=\bt\omega^2_t - \bt d\omega^1
=\partial_t\bt\omega^2 - d\bt\omega^1 =0,
\]
i.e., the  Neumann data for the
function $*\omega^3$ vanish at $\p M$. Thus, $\omega^3(t)=0$.
\hfill$\Box$

Assume that $\omega(t)$ is a smooth solution of the
complete system (\ref{ibvp}).
The {\em complete Cauchy data} of $\omega(t)$ consist of
\[
  \big(\bt\omega(x,t),\bn\omega(x,t)\big)
  ,\quad (x,t)\in\partial M\times\R.
\]
The Cauchy data for the solutions $\omega(t)$
of Maxwell's equations have a particular
structure. Indeed, by taking the tangential trace of the equation (\ref{MF1}),
we obtain  $\bt\omega^2_t= -d\bt\omega^1$. Further, by integrating,
\beq\label{27-}
  \bt\omega^2(x,t)=\bt \omega^2(x,0)-\int_0^t d(\bt\omega^1(x,t'))\,dt',
\quad x\in\partial M.
\eeq
Similarly, by taking the normal trace of
equation (\ref{MA1}), we find that $\bn \omega_t^1=d\bn \omega^2$,
so by integrating
\beq
\label{2.12.1}
\bn\omega^1(x,t)= \bn \omega^1(x,0)+
\int_0^t d(\bn\omega^2(x,t'))\,dt',
\quad x\in\partial M.
\eeq
In this work, we consider mainly the case $\omega(0)=0$, when the
lateral Cauchy data for the original problem
of electrodynamics is simply
\begin{eqnarray}
\bt\omega &=& (0,f,-\int_0^td f(t')dt'),
\label{t-data}\\
\bn\omega &=& (0,g,\int_0^td g(t')dt'),
\label{n-data}
\end{eqnarray}
where  $f$ and
$g$ are functions of $t$ with values in $\Omega^1\partial M$.
The following
theorem implies that solutions of Maxwell's equations
are solutions of the complete Maxwell system and gives sufficient
conditions for the converse result.

\begin{theorem}\label{equivalence}
If $\omega(t)\in C(\R,{\bf H})\cap
C^1(\R,{\bf L}^2)$
satisfies the equation
\begin{equation}\label{M-eq}
  \omega_t +{\mathcal M}\omega = 0,\quad t>0,
\end{equation}
with $\omega(0)=0$,
and $\omega^0(t)=0$, $\omega^3(t)=0$, then $\bt\omega,\, \bn\omega$
are of the form (\ref{t-data})--(\ref{n-data}).

Conversely, if $\bt\omega,\, \bn\omega$ are of  the form
(\ref{t-data})--(\ref{n-data})
for $0 \leq t \leq T$,
and $\omega(t)$ satisfies
  (\ref{M-eq}), with $\omega(0)=0$,
  then $\omega(t)$ is a solution of Maxwell's equations, i.e.,
$\omega^0(t)=0$, $\omega^3(t)=0$.
\end{theorem}

{\em Proof:}
The first part of the theorem follows
from the above considerations if
  we show that $\omega(t)$ is sufficiently regular.
For  $\omega^2 \in C( \R,H(\delta, \Omega^2M))$, by Proposition \ref{paquet},
${\bf n}\omega^2 \in C( \R,H^{-1/2}(\Omega^1 \partial M))$ with
$d{\bf n}\omega^2 \in C( \R,H^{-1/2}(\Omega^2 \partial M))$. As
$\delta \omega_t^1(t) = \delta \delta \omega^2(t) =0$, it holds also that
\[
{\bf n}\omega_t^1 \in C( \R,H^{-1/2}(\Omega^2 \partial M)), \quad
\]
implying (\ref{n-data}).
To prove (\ref{t-data}), we use Maxwell duality: Consider the forms
  \[
  \eta^{3-k}=(-1)^k*\alpha^{-1}\omega^k.
\]
Then $\eta=(\eta^0,\eta^1,\eta^2,\eta^3)$
satisfies the complete system dual to the Maxwell system,
$\eta_t+\tilde \M \eta=0$, where
$\tilde \M$ is the Maxwell operator with metric $g$
and scalar impedance $\alpha^{-1}$.
Then formula (\ref{n-data}) for the
solution $\eta$ implies  (\ref{t-data}) for $\omega$.

To prove the converse, we observe that the equations,
\begin{eqnarray}\label{aux}
\p_t \omega^0(t) -\delta\omega^1(t) = 0,\quad
\p_t \omega^1(t) +d\omega^0(t) -\delta\omega^2(t)
=0
\end{eqnarray}
imply that
\[
  \omega_{tt}^0(t) +\delta d\omega^0(t)=0.
\]
In addition, $\omega^0(0)=0$, $\omega^0_t(0)=0$,
and from (\ref{t-data}),  ${\bf t}\omega^0(t) =0$. Thus,
$\omega^0 = 0$ for $0\leq t \leq T$.
By the Maxwell duality
described earlier, this implies also that $\omega^3(t)=0$.
\hfill$\Box$

The following definition, where $R$ is a right inverse to the
mapping $\bt$,
fixes the solutions of
the forward problem used in this work.

\begin{definition}\label{def. weak Maxwell}
Let $h=(h^0,h^1,h^2)\in C^{\infty}([0,T], {\bf \Omega}\partial M)$.
The solution $\omega(t)$ of the initial boundary
value problem
\ba
& &\omega_t+ {\mathcal M}\omega=0,\quad t>0, \\
\noalign{\vskip4pt}
& & \omega(0)=\omega_0\in {\bf L}^2(M),\quad \bt\omega = h,
\ea
is given by
\[
  \omega(t) = Rh(t) + {\mathcal U}(t)\omega_0
   - \int_0^t{\mathcal U}(t-s)({\mathcal M}Rh(s)
+Rh_s(s))ds.
\]

When $\omega_0=0$ and $h$
is a smooth boundary source of the form
(\ref{t-data}),
\ba
h= (0,f,-\int_0^t df(t')dt'),\quad f\in C^\infty_0(]0,T[,\Omega^1\p M),
\ea
$\omega(t)$ is called the  solution of
Maxwell's equations in $M\times [0,T]$ with the boundary condition
$\bt \omega^1=f$ and the initial condition $\omega(0)=0$.
\end{definition}

To emphasize the dependence of  $\omega(t)$ on $f$ above, we  write
occasionally
\beq\label{formul 1}
\omega(t)=\omega^f(t) = (0,(\omega^f)^1,(\omega^f)^2,0).
\eeq

We note that  $f$ could be
chosen from a wider class, e.g. from
$H^{1/2}(\partial M\times [0,T])$.

We use the notation
$\Cnull([0,T],\Omega^1\partial M)$ for the space of $C^{\infty}$
functions $[0,T]\to\Omega^1\partial M$ vanishing near $t=0$.
Theorem \ref{equivalence} motivates the following definition.

\begin{definition}
\label{27.11.d}
The admittance map, ${\mathcal Z}^T$ is defined as
\[
{\mathcal Z}^T: \Cnull([0,T],\Omega^1\partial M)\to
  \Cnull([0,T],\Omega^1\partial M),\quad
\bt \omega^1|_{\partial M \times [0,T]} \mapsto\bn
\omega^2|_{\partial M \times [0,T]},
\]
where $\omega(t)$ is the solution of
Maxwell's equations in $M\times [0,T]$ with
$\omega(0)=0$.
\end{definition}

Note that in the
classical terminology of electric and magnetic fields,
${\mathcal Z}^T$ maps the tangential electric field
$n\times E|_{\partial M\times [0,T]}$ to the tangential magnetic field
$n\times H|_{\partial M\times [0,T]}$.

The following result, which relates the boundary data and the energy of the
electromagnetic field, is
crucial  for  boundary control. It is
a
version of the Blagovestchenskii
formula (see \cite{BeBl} for the scalar case).

\begin{theorem}\label{blacho}
\begin{enumerate}

\item  For any $T >0$ and $f,\, h \in \Cnull([0,2T],\Omega^1\partial M)$,
the knowledge of the admittance map $\cZ^{2T}$ allows us to
evaluate the inner products
\[
  ((\omega^f)^j(t),(\omega^h)^j(s))_{L^2},\quad j=1,2,\quad
  0\leq s,t\leq T.
\]

\item
For any $T >0$ and  $f\in\Cnull([0,T],\Omega^1(\partial M))$,
$\cZ^{T}$ determines the energy of the field $\omega^f$ at $t=T$,
defined as
\[
{\mathcal E}^{T}(\omega^f)= \frac 12 \|(\omega^f)^1(T)\|^2_{L^2}+
\frac 12 \|(\omega^f)^2(T)\|^2_{L^2}.
\]
\end{enumerate}
\end{theorem}

{\em Proof:} {\em 1.}
Let $\omega(t)=\omega^f(t)$ and $\eta(s)=\omega^h(s)$
with
$
  F^j(s,t) = (\omega^j(s),\eta^j(t))_{L^2},$ $ j=1,2.
$
By (\ref{complete}),
\begin{eqnarray}\label{A-fomrmu}
& & (\partial_s^2 -\partial_t^2)F^j(s,t) =
   (\omega_{ss}^j(s),\eta^j(t))_{L^2}-
   (\omega^j(s),\eta_{tt}^j(t))_{L^2}\\
  \nonumber
&=&-((d\delta+\delta d)\omega^j(s),\eta^j(t))_{L^2}+
  (\omega^j(s),(d\delta+\delta d)\eta^j(t))_{L^2} =
  b^j(s,t).
\end{eqnarray}
We apply Maxwell's equations (\ref{MF}), (\ref{MA}) and the
commutation relations,
\beq \label{commut}
\bt d\omega^j = d\bt\omega^j,\quad
  \bn\delta\omega^j = \bt **d\frac 1\alpha*\omega^j =
  d\bt\frac 1\alpha*\omega^j = d\bn \omega^j,\quad j=1,2,
\eeq
where $d$, in the right side is the exterior derivative on $\p M$.
A straightforward applications of Stokes' formula (\ref{stokes1}), yields
\begin{eqnarray*}
   b^1(s,t)= \langle \bn\omega^2_s(s),\bt\eta^1(t)\rangle
- \langle \bt\omega^1(s),\bn\eta^2_t(t)\rangle,\\
  b^2(s,t) = \langle \bn\omega^2(s),\bt\eta^1_t(t)\rangle
- \langle \bt\omega^1_s(s),\bn\eta^2(t)\rangle.
\end{eqnarray*}
As $\cZ^{2T}$ determines $b^1(s,t)$ and $b^2(s,t)$ for $t,s<2T$ and
\beq\label{B-fo}
  F^j(0,t)=F^j(s,0)=0,\quad F_s^j(0,t)=F_t^j(s,0)=0,
\eeq
the function  $F^j(s,t)$ can be found from
the wave equation (\ref{A-fomrmu}) for $s+t <2T$.

{\em 2.} Again, by differentiating and using Maxwell's equations
and Stokes' formula, we obtain
\bfo
\p_t{\mathcal E}^t(\omega^f) = -\langle\bt\omega^1(t),\bn\omega^2(t)\rangle
=  -\langle f(t),\cZ^Tf(t)\rangle.
\efo
As ${\mathcal E}^0(\omega^f)=0$, the energy is readily obtained for
$t\leq T$.
\hfill$\Box$

\subsection{Unique continuation results}\label{sect: continuation}

For further applications to inverse problems, in this section we
consider the unique continuation of the Holmgren-John type
for Maxwell's equations. We start with an extension of
differential forms outside the manifold $M$. Let $\Gamma\subset\partial M$
be open and $\tilde M$ be
an extension of $M$ across $\Gamma$, i.e.,
$M\subset\tilde M$, $\Gamma\subset \tilde M^{\rm int}$
and $\partial M\setminus\Gamma\subset\partial
\tilde M$. Let  $\tilde{g}, \, \tilde{\alpha}$ be smooth
continuations of   $g$ and $\alpha$ to $\tilde M$.
In this case, we say that the manifold $(\tilde M, \tilde{g},\tilde{\alpha})$
is an {\em extension of $(M,g,\alpha)$ across
$\Gamma$.}

Let $\omega^k$ be a $k$-form  on $M$ and $\tilde\omega^k$ its extension
by zero  to $\tilde M$. It follows from Stokes' formula (\ref{stokes1})
that for  $\omega^k\in H(d,\Omega^k M)$ with
$\bt\omega^k\big|_\Gamma=0$, we have $\tilde \omega^k\in
H(d,\Omega^k\tilde M)$. Similarly, if $\omega^k\in H(\delta,
\Omega^k M)$ and
$\bn\omega^k\big|_\Gamma=0$, then $\tilde
\omega^k\in H(\delta,\Omega^k\tilde M)$.
These yield the
following result.

\begin{proposition} \label{proposition1}
Let $\omega(t)\in C^1(\R,{\bf L}^2)\cap C(\R,{\bf H})$
with $\bt\omega|_{\Gamma \times [0,T]}=0$ and
$\bn\omega|_{\Gamma \times [0,T]}=0$, be a solution of
the complete Maxwell system (\ref{complete})
in $M\times [0,T]$. Let $\tilde\omega$
be its extension by zero across $\Gamma\subset  \partial M$. Then
$\tilde{\omega}(t)$ satisfies the complete Maxwell system (\ref{complete}) in
$\tilde M\times [0,T]$.
\end{proposition}

We are particularly interested in the solutions
of Maxwell's equations. The following
result extends  Proposition \ref{proposition1}  to this case.

\begin{lemma}\label{2.12l}
Let $\omega(t) \in C^1(\R,{\bf L}^2)\cap C(\R,{\bf H})$
be a solution of Maxwell's equations (\ref{MF}), (\ref{MA})
in $M\times [0,T]$, i.e.,
$\omega^0(t)=0, \,\omega^3(t)=0$. In addition, let
$\bt\omega^1|_{\Gamma \times [0,T]}=0$,
$\bn\omega^2|_{\Gamma \times [0,T]}=0$, and $\omega(0)=0$. Then
$\tilde{\omega}(t) \in C^1(\R,{\bf L}^2(\tilde{M}))\cap
C(\R,{\bf H}(\tilde{M}))$
is a solution of Maxwell's equations in $\tilde M\times [0,T]$.
\end{lemma}

{\em Proof:} The above conditions together with
Theorem \ref{equivalence} imply  that
\[
  \bt \omega =
(0,\bt\omega^1,-\int_0^td\bt\omega^1dt') =0,\quad
  \bn \omega  =
  (0,\bn\omega^2,\int_0^td\bn\omega^2 dt') =0
\]
in $\Gamma\times[0,T]$.
Therefore, by  Proposition \ref{proposition1},
  $\tilde \omega (t)$ satisfies   (\ref{complete}) in
$\tilde M\times [0,T]$.
Clearly, also $\tilde \omega^0(t)=0$,  $\tilde \omega^3(t)=0$ in
$\tilde M\times [0,T],$ and $\tilde \omega(0)=0$.

\hfill$\Box$

When we deal with a general solution to  Maxwell's equations
(\ref{MF})--(\ref{MA}), which may not satisfy zero initial conditions,
and try to extend them by zero across $\Gamma$, the arguments of
Lemma  \ref{2.12l} fail.  Indeed, if $\omega(0) \neq 0$, then
(\ref{2.12.1}) show that $\bn \omega^2 =0$ is not sufficient
for $\bn \omega^1 =0$. However, by
differentiating with respect to $t$, the
parasite term $\bn\omega^1(0)$ vanishes.
This is the motivation why
Theorem \ref{UCP}  below deals with the time
derivatives of weak solutions.

Let, again,  $\Gamma\subset\partial M$ be open and $T >0$.
Denote by $K(\Gamma,T)$ the double cone of influence with the base
on the slice $t=T$,
\[
  K(\Gamma,T)=\{(x,t)\in M\times [0,2T]\mid
  \tau(x,\Gamma)<T-|T-t|\}
\]
where $\tau(x,y)$ is the distance function on $(M,g)$ (see Figure 1).

\begin{figure}[htbp]
\begin{center}
\psfrag{1}{$\Gamma$}
\includegraphics[width=8cm]{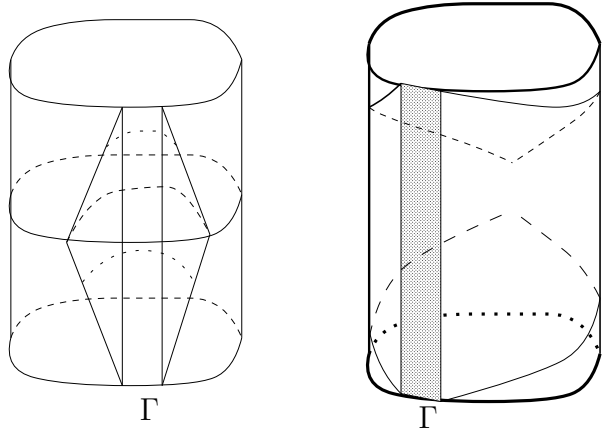} \label{pic 2}
\end{center}
\caption{Left: The double cone of influence. Right:
For $T$ large enough, the double cone contains a slice $\{T/2\}\times M$.}
\end{figure}

We prove the following unique continuation result for the time derivative
of the fields.

\begin{theorem}\label{UCP}
Let $\omega(t)$ be a weak solution (\ref{weak solution}) of
the initial boundary value problem (\ref{ibvp}).
Assume that $\omega_0 =(0,\omega_0^1,
\omega^2_0,0)$, $\delta\omega^1_0=0$,
$d\omega_0^2=0,$ and $\rho =0$.
If $\bn\omega^2=0$ on $\Gamma\times]0,2T[$, then
$\omega_t=0$ in  $K(\Gamma,T)$.
\end{theorem}

{\em Proof:}
When $\omega(t) \in C^2(]0, 2T[, {\bf L}^2) \cap C^1(]0, 2T[, {\bf H})$,
then $\eta(t) = \omega_t(t) \in C^1(]0, 2T[, {\bf L}^2) \cap C^0(]0, 2T[,
{\bf H})$ also satisfies Maxwell's equations (\ref{MF}), (\ref{MA}).
{\newtext Let $\tilde M$ be
the extension of $M$ across
$\Gamma$ and $\tilde \eta$ be the  extension of $\eta$ by zero.
In follows for (\ref{27-}) and (\ref{2.12.1}) that
$\bt \eta^2 =-d \bt \omega^1 =0$ and
$\bn \eta^1 = d \bn \omega^2 =0$ in $\Gamma\times]0,2T[$.
Therefore, by Proposition \ref{proposition1},
\[
\tilde{\eta}(t)  \in C^1(]0, 2T[, {\bf L}^2(\tilde{M}))
\cap C^0(]0, 2T[, {\bf H}(\tilde{M})),
\]
is a solution of the complete system and obviously, also a
solution of Maxwell's equations (\ref{MF}), (\ref{MA})
in $\tilde{M} \times]0,2T[$. }

By the unique continuation for sufficiently smooth solutions
(see \cite{EIsNkTa} and Remark \ref{remark Tataru} below), we have that, for any $\sigma >0$,  $\tilde\eta=0$
in the double cone
\[
\{ (x,t), \, x \in \tilde M, t \in\R| \,  \tilde{\tau}(x,\tilde M
\setminus M)<T - \sigma-|T-t| \}.
\]
  Thus, $\eta = 0$ in $K(\Gamma,T)$.

When $\omega(t) \in C^0(]0, 2T[, {\bf L}^2)$ is a weak solution,
we use  Friedrich's mollifier in $t$,
\bfo
\omega_{\sigma} = \psi_{\sigma} * \omega, \quad
\psi_{\sigma}(t) =  (1/\sigma) \psi(t/\sigma) \quad
\hbox{for $\sigma >0$},
\efo
where $\psi \in C^{\infty}_0([-1,1])$, $\int \psi(s)ds=1$. Then
$\omega_{\sigma}(t)$ satisfies
conditions of the theorem with $\Gamma\times]0,2T[$ replaced by
$ \Gamma \times ]\sigma, 2T-\sigma[.$ As
\[
{\mathcal M}^j \omega_{\sigma}  = (-\partial _t )^j\omega_{\sigma}
  \in C^{\infty}(]\sigma, 2T-\sigma[, {\bf L}^2(M)), \quad
\hbox{for any $j>0$},
\]
we have in particular that $\omega_{\sigma}(t) \in C^2(]\sigma,
2T-\sigma[, {\bf L}^2)
\cap C^1(]\sigma, 2T-\sigma[, {\bf H})$. By the above,
$\p_t\omega_{\sigma}(t) =0$
in $K_{\sigma}(\Gamma,T)$, where
\bfo
K_{\sigma}(\Gamma,T) = \{(x,t)\mid \tilde{\tau}(x,\tilde{M}
\setminus M) < T-\sigma-|T-t|\}.
\efo
As $\p_t\omega_{\sigma}(t) \to \p_t\omega(t)$ in the distribution sense,
the result follows.
\hfill$\Box$

\begin{remark}\label{remark Tataru}
{\rm The article \cite{EIsNkTa}, based on results of
Tataru  \cite{Ta1,Ta3} deals with
scalar $\epsilon$ and $\mu$. However, due to the single velocity of the wave
propagation, it is, in principle, possible to generalize it to the scalar impedance case.
Another way to prove the desired unique continuation for sufficiently 
smooth solutions of the equation (\ref{wave equation}) is to use
the simplified version of Tataru's construction, given in 
\cite[sec. 2.5]{KKL}. 
There, the unique continuation result 
is based on local Carleman estimates for the solutions of the scalar wave equation,
$u_{tt}- a_{ij}(x)\p_{i}\p_{j}u+A_1(x,D)u=0$, where $A_1(x,D)$ is the first-order differential operator.
These estimates utilized a function $\phi(x,t)$
that is pseudoconvex with respect to the metric $a_{ij}$,
and absorbed the perturbation due to $A_1(x,D)$ into the main terms 
of the Carleman
estimates. By  (\ref{Laplace}) and (\ref{wave equation}),
  the operator
${\mathcal M}^2$ is in local coordinates a principally diagonal operator 
with the same second order differential operator,
$g^{ij}\p_{i}\p_{j}$, acting on all components of $\omega(t)$. 
As in \cite{KKL}, one can treat the first-order terms
as a perturbation and obtain the desired Carleman
estimate. 
In this manner, the constructions in \cite{KKL} can be word-by-word generalized to the considered
case of Maxwell's equations with scalar wave impedance.
}
\end{remark}

\begin{remark}\label{remark 1a}
{\rm
It is clear from the above arguments that if $\omega(t)$ is a weak solution of
the initial boundary value problem (\ref{ibvp}) and
$\omega(t) \in C^{\infty}(]0, 2T[, {\bf L}^2(M))$,
then
\[
\omega(t)  \in C^{\infty}(]0, 2T[, {\mathcal D}^{\infty}
({\mathcal M}_{\rm e})), \quad
\omega(t)  \in C^{\infty}(M^{\rm int} \times ]0, 2T[),
\]
where we used the notation
${\mathcal D}^{\infty}({\mathcal M}_{{\rm e}}) = \bigcap_{N>0}
  {\mathcal D}({\mathcal M}^N_{{\rm e}})$.}
\end{remark}

\subsection{Controllability results}

In this section we derive controllability
results for Maxwell's equations. We divide these results
in {\em local results}, i.e., controllability at short times and
{\em global results},
where the time of control is long enough so that
the controlled  electromagnetic
waves fill the whole manifold. Both types of
results are based
on the unique continuation of Theorem \ref{UCP}.

Let $\omega^f(t)$, $f\in C^\infty_0
(\R_+, \Omega^1\partial M)$ be a solution of Maxwell's equations
in the sense of Definition \ref{def. weak Maxwell}
with the initial condition $\omega^f(0)=0$.
Let $\tilde\omega$ be the weak solution of (\ref{ibvp}) given by
(\ref{weak solution}) with $\rho=0$ and $\tilde \omega(T)=\omega_0=
(0,\omega_0^1,\omega_0^2,0)$. Similar considerations to
those in the proof of Theorem \ref{blacho},
show that
\begin{equation}\label{control identity}
(\omega^f(T),\omega_0)_{\bf L^2}
=-\int_0^T\langle\bt\omega^f(t),\bn\tilde\omega(t)
\rangle dt,
\end{equation}
which we will refer to as the {\em control identity}.

\subsubsection{Local controllability}\label{local controllability
section}

In this section, we study differential $1-$forms in $M$
generated by boundary sources active for short periods of time.
Instead of a complete characterization of these forms,
we show that they form a sufficiently large  subspace
in $ L^2(\Omega^1M)$. The difficulty that prevents a complete
characterization of this subspace lies in  topology of
the domains of influence, which can be very complicated.

Let $\Gamma\subset\partial M$ be open and $T>0$. The {\em domain of
influence}, $ M(\Gamma,T)$, is defined as
\beq
\label{domain of influence}
  M(\Gamma,T)=\{x\in M| \tau(x,\Gamma)< T\}, \\
\nonumber
  M(\Gamma,T)\hspace{-1mm}\times\hspace{-1mm}
\{T\}= K(\Gamma,T) \cap \{t=T\}.
\eeq

Let $C^\infty_0(]0,T[,\Omega^1\Gamma) \subset
C^\infty_0(]0,T[,\Omega^1\partial M)$ consists of forms supported
in $\overline{\Gamma}\times [0,T]$ and
\beq\label{3.12.2}
  X(\Gamma,T)={\rm cl}_{L^2}\{ (\omega^f)^1(T)\mid f\in
  C^\infty_0(]0,T[, \Omega^1\Gamma )\}.
\eeq
Furthermore, let
\[
H(\delta,M(\Gamma,T)) =\{
\omega^2\in H( \delta,\Omega^2M)\mid {\rm supp}\,(\omega^2)
\in \overline{M(\Gamma,T)}\}.
\]
For $ S \subset M$, we define
  $\,H^1_0(\Omega^kS) \subset H^1_0(\Omega^kM)$ consisting of
$k$--forms with support in $\overline{S}$.

\begin{theorem}\label{local control th}
For any open $\Gamma\subset\partial M$  and $T>0$,
\beq
\label{1a2}
  \delta H^1_0(\Omega^2M(\Gamma,T))\subset
  X(\Gamma,T)\subset{\rm cl}_{L^2}\left(\delta H(
\delta,M(\Gamma,T)
)\right).
\eeq
\end{theorem}

{\em Proof:}
The right inclusion being an immediate corollary of (\ref{MF}),
we concentrate on the left one.

Let $\omega_0^1\in L^2(\Omega^1M)$
satisfy
\beq\label{1a3}
  \left(\omega_0^1,(\omega^f)^1(T)\right)_{L^2}=0,
\quad \hbox{for all } f\in C^\infty_0(]0,T[,\Omega^1\Gamma).
\eeq
By the Hodge decomposition (see \cite{Sc}) in $L^2(\Omega^1M)$, we have
\beq
\label{1a4}
\omega^1_0=\hat \omega^1_0+\delta \eta^2_0,
\eeq
where $d\hat \omega^1_0=0$, $\bt \hat \omega^1_0=0$ and $\eta^2_0\in
H(\delta,\Omega^2 M)$. Thus, (\ref{1a3}) is equivalent to
\beq
\label{1a5}
  \left(\delta \eta^2_0,(\omega^f)^1(T)\right)_{L^2}=0.
\eeq
Let $\tilde\omega(t)$ be a weak solution to (\ref{ibvp}) with $\rho=0$
and the initial data at $t=T$ given by
$ \tilde\omega(T)= (0,\delta \eta^2_0,0,0)$.
By the control identity (\ref{control identity}), the
orthogonality  (\ref{1a3}) and the particular form of the
boundary data for solutions of Maxwell's equations
(\ref{t-data}), (\ref{n-data}), we see that
\[
0= \int_0^T\langle \bt \omega^f(t),\bn\tilde\omega(t)\rangle
  =\int_0^T\langle \bt(\omega^f)^1(t),\bn\tilde\omega^2(t)
\rangle = \int_0^T\langle f,\bn\tilde\omega^2(t)
\rangle,
\]
i.e., $\bn\tilde\omega^2=0$ on $\Gamma\times ]0,T[$.  Since
\[
  \tilde\omega(T+t)=(0,\tilde\omega^1(T-t),
  -\tilde\omega^2(T-t),0),
\]
also  $ \bn\tilde\omega^2=0$ on $\Gamma
\times ]T,2T[$. Since $\delta \tilde \omega^2(t)=0$,
we see by using Proposition \ref{paquet} that  $\bn\tilde\omega^2 \in
C^0(\R,H^{-1/2}(\Omega^1\p M))$. Hence
\[
\bn\tilde\omega^2=0 \quad \hbox{on} \,\,\,\Gamma
\times ]0,2T[.
\]
Therefore, by Theorem \ref{UCP}, $\p_t\tilde\omega^2=0$ in $K(\Gamma,T)$,
In particular,
$d\delta \eta^2_0=-\p_t \tilde\omega^2(T)=0$ in $M(\Gamma,T)$.
In other words, if
$\omega^1_0 \in X(\Gamma,T)^\perp$, then
the term $\eta^2_0$ in the decomposition (\ref{1a4}) satisfies
$d\delta \eta^2_0=0$ in $M(\Gamma,T)$.
For any $\nu^2 \in H^1_0(\Omega^2M(\Gamma,T))$, we have therefore
\[
\left(\delta \nu^2,\omega_0^1\right)_{L^2}=
\left(\nu^2,d\delta \eta^2_0\right)_{L^2}=0,
\]
and thus $\delta \nu^2 \in(X(\Gamma,T)^\perp)^\perp$.
This is equivalent to the left inclusion in
(\ref{1a2}).
\hfill$\Box$

\begin{remark}
\label{Remark 1.}{\rm Later in this work, {\newtext we deal mainly with
the time derivatives of
the electromagnetic  fields. Since
$\omega_t^f(t)=\omega^{\p_t f}(t)$, we see by using
\beq\label{a fin}
  X(\Gamma, T) \subset {\rm cl}_{L^2}\{(\omega^f_t(T))^1\mid f\in \Cnull
([0,T],\Omega^1\p M)\}
\eeq
and (\ref{MF}), that the inclusions (\ref{1a2}) remain valid when
$ X(\Gamma, T)$ is replaced with the right hand side of (\ref{a fin}).}}
\end{remark}

\subsubsection{Global controllability}

This section is devoted to the study of controllability results
when the control times are large enough so that the waves fill the whole
manifold.

For $\Gamma \subset \p M$ and $T>0$, we define
\beq
\label{3.12.3}
  Y(\Gamma,T) = \{\omega^f_t(T)\mid
  f\in C^\infty_0(]0,T[, \Omega^1\Gamma)\},
\eeq 
 where
$\Omega^1\Gamma$ is the set of 1-forms in $\Omega^1\p 
M$
supported on $\Gamma$
and abbreviate  $Y(\partial M,T)=Y(T)$. Our objective
is to characterize  $Y(\Gamma,T)$ for $T$ large enough.
In the following theorem, we use the notation
\beq
\label{23.6.1}
{\rm rad}_\Gamma\,(M)=\max_{x\in M} \tau(x,\Gamma),\quad
{\rm 
rad}\,(M)={\rm rad}_{\p M}\,(M).
\eeq

\begin{theorem}\label{global control th}
For open non-empty $\Gamma\subset \p M$ and $T\geq T_0> 2\,{\rm 
rad}_\Gamma(M)$, we have
${\rm cl}_{{\bf L}^2(M)}Y(\Gamma,T)=Y$, where
\beq
\label{3.12.5}
Y = \{0\}\times\delta H(\delta)\times
  d\Hnull(d)\times\{0\}.
\eeq
\end{theorem}

\noindent
{\em Proof:} Let $\omega(t)=\omega^f(t)$ be a solution,
in the sense of Definition
\ref{def. weak Maxwell}, of the initial boundary-value problem
  with  $f\in C^\infty_{0}
([0,T_0],\Omega^1\Gamma)$. Since $f=0$ for $T\geq T_0$,
we have $\bt\omega^1(T)=0$,
and consequently,
\begin{eqnarray*}
  \omega_t(T) = -{\mathcal M}\omega(T)
  = (0,\delta\omega^2(T),-d\omega^1(T),0) \in\{0\}\times\delta H(\delta)\times
  d\Hnull (d)\times\{0\}.
\end{eqnarray*}
To prove the converse, we will show that
$Y(\Gamma,T)$ is dense in $\{0\}\times\delta H(\delta)\times
  d\Hnull (d)\times\{0\}$. To this end,
let $\omega_0 \in\{0\}\times\delta H(\delta)\times
  d\Hnull (d)\times\{0\}$ and
   $\omega_0\perp Y(\Gamma,T)$, i.e.,
\begin{equation}\label{orthogonality}
(\omega_0,\omega_t(T))_{{\bf L}^2} =(\omega_0^1,\omega_t^1(T))_{L^2}
+(\omega_0^2,\omega_t^2(T))_{L^2} =0,
\end{equation}
for any $f\in C^\infty_{0}
([0,T_0],\Omega^1\Gamma )$.

Let $\tilde{\omega}$ be the weak solution of the problem
\ba
  \tilde{\omega}_t + {\mathcal M}\tilde{\omega} =0, \quad
  \bt \tilde{\omega} =0,\quad \tilde{\omega}(T)=\omega_0.
\ea
Observe that  $\omega_0$ satisfies
$
  \delta\omega_0^1 =0$ and $d\omega_0^2 =0,
$
so that
  $\tilde{\omega}$ satisfies Maxwell's equations.
Consider the function $F:\R\to\R$,
$
   F(t) = (\tilde{\omega}(t),\omega_t(t))_{{\bf L}^2}.
$
By  Maxwell's equations,
\begin{eqnarray*}
  F_t(t)&=& (\tilde{\omega}(t),\omega_{tt}(t))_{{\bf L}^2} +
(\tilde{\omega}_t(t),
\omega_t(t))_{{\bf L}^2} \\
  &=& -(\tilde{\omega}^1,\delta d\omega^1)_{L^2} - (\tilde{\omega}^2,
  d\delta\omega^2)_{L^2} +(d\tilde{\omega}^1,d\omega^1)_{L^2} + (\delta
\tilde{\omega}^2,\delta\omega^2)_{L^2},
\end{eqnarray*}
and further, by Stokes' formula (\ref{stokes1}),
\[
  F_t(t) = -\langle\bt \tilde{\omega}^1(t),\bn d\omega^1(t)\rangle
  -\langle\bn \tilde{\omega}^2(t),\bt\delta\omega^2(t)\rangle.
\]
However, ${\bt}\tilde \omega =0$ and $\delta \omega^2 =\omega^1_t$. Thus,
\[
  F_t(t) = -\langle\bn \tilde{\omega}^2(t),\bt\omega^1_t(t)\rangle
  =-\langle\bn \tilde{\omega}^2(t),f_t(t)\rangle.
\]
On the other hand, since $\omega(0)=0$,
  the orthogonality condition
(\ref{orthogonality}) implies
that $F(0)=F(T)=0$,  i.e.,
\[
  \int_0^T\langle\bn \tilde{\omega}^2(t),f_t(t)\rangle dt = -\int_0^T
  F_t(t)dt =0.
\]
Since
  $\, f \in C^{\infty}_0(]0,T[, \Omega^1 \Gamma)$ is arbitrary, this implies
  that
$\bn \tilde{\omega}^2_t =0$ in $\Gamma\times]0,T[$.
Thus, by Theorem \ref{UCP},
$\tilde{\omega}_{tt}= 0$ in the double cone $ K(\Gamma, T/2)$.
Since  $T_0> 2\,{\rm rad}_\Gamma(M)$, this double cone contains
a cylinder  of the form $C = M\times]T/2-s,T/2+s[$ with some $s>0$.
(See Figure 1).

As $\tilde{\omega}_{tt}$ satisfies Maxwell's equations with
a homogeneous boundary condition $\bt \tilde{\omega}_{tt}=0$,
this implies that $\tilde{\omega}_{tt}= 0$ in $M \times \R.$ Therefore,
$
  \tilde{\omega}(t)=\omega_1 + t\omega_2,
$
where $\omega_1$ and $\omega_2$ do not depend on $t$.
Again, by Maxwell's equations,
\[
  \omega_2 =\tilde \omega_t ={\mathcal M}\omega_1 +
  t{\mathcal M}\omega_2.
\]
Therefore,
$
  \omega_2 ={\mathcal M}\omega_1$ and ${\mathcal M}\omega_2 =0.
$
But then Stokes' formula implies that
\[
  (\omega_2,\omega_2)_{{\bf L}^2} = (\omega_2,{\mathcal M}\omega_1)_{{\bf L}^2}
  = -({\mathcal M}\omega_2,\omega_1)_{{\bf L}^2}=0,
\]
i.e., $\omega_2=0$ and  ${\mathcal M}\omega_1=0$.
Furthermore, by the choice of $\omega_0$,
\[
  \omega_1 = \tilde{\omega}(T)=\omega_0 = (0,-\delta\nu^2,
d\nu^1,0)= {\mathcal M}\nu,
\]
with $\nu\in \{0\}\times \Hnull (d)\times H(\delta)
\times\{0\}$. By a further application of
Stokes' formula,
\[
(\omega_1,\omega_1)_{{\bf L}^2} = (\omega_1,{\mathcal M}\nu)_{{\bf L}^2}
  = -({\mathcal M}\omega_1,\nu)_{{\bf L}^2}=0,
\]
i.e., $\omega_1=0$ and, therefore, $\omega_0=0$.
\hfill$\Box$

\subsection{Generalized sources}\label{generalized sources section}

So far, we dealt only with smooth boundary
sources and the corresponding fields. Later, we need more general
sources which are described in this section.

Let $W^T$ be the wave operator,
\[
  W^T:C^\infty_0(]0,T_0[, \Omega^1\partial M )\to Y,
\quad f\mapsto \omega_t^f(T),
\]
where $T\geq T_0>2\,{\rm rad}\,(M)$.
Let $\|\cdot\|_{{\mathcal F}}$ be a quasinorm on the
space of boundary sources defined via $W^T$,
\begin{equation}\label{f norm}
  \|f\|_{{\mathcal F}} = \|W^T f\|_{{\bf L}^2}.
\end{equation}
  By energy conservation,  this norm is independent of $T\geq T_0$ and
by Theorem \ref{blacho},
if the admittance map ${\mathcal Z}^T$ is given, we can
evaluate $ \|f\|_{{\mathcal F}}$ for $f \in C^\infty_0(]0,T_0[,
\Omega^1\partial M )$.

Using the standard procedure in PDE-control, e.g.  \cite{Ru,LTr},
there is a Hilbert space of {\em generalized boundary sources} with the
norm defined by (\ref{f norm}). Indeed, we first introduce the
space $\F([0,T_0])$,
\bfo
\F([0,T_0])=C^\infty_0(]0,T_0[, \Omega^1\partial M )/\sim,
\efo
where $
  f\sim g\mbox{ iff $W^Tf=W^Tg$},
$
and then complete it with respect to the norm (\ref{f norm})
to obtain ${\overline{\mathcal F}}([0,T_0])$.
By Theorem \ref{global control th}, $W^T$ is an isometry
between ${\overline{\mathcal F}}([0,T_0])$
and $Y$ for any $T\geq T_0>2\,{\rm rad}\,(M)$.
{\newtext The elements of ${\overline{\mathcal F}}([0,T_0])$
are equivalence classes of Cauchy sequences $(f_j)_{j=1}^\infty$
and we denote them by $\hat{f}=(f_j)_{j=1}^\infty$. (This is
  slight abuse of notations, as $(f_j)_{j=1}^\infty$
is a representative of the equivalence class $\hat{f}$.)}
To put it in another
way, for any
$\omega_0 \in Y$, there is a sequence $(f_j)_{j=0}^\infty$ with $f_j\in
C^\infty_0(]0,T_0[, \Omega^1\partial M )$, defining a generalized source
$\hat{f} \in {\overline{\mathcal F}}([0,T_0])$, and for the corresponding
wave
\beq
\label{gener}
  \omega_t^{\hat f}(t) := \lim_{j\to \infty}\omega_t^{f_j}(t), \quad
\hbox{for }t\geq T_0,
\eeq
we have
$
\omega_t^{\hat f}(T)=\omega_0.
$
Since in this work $T_0$ is considered as a fixed parameter, we
denote  ${\overline{\mathcal F}}([0,T_0])$ for brevity as
$\overline{\F}$.

We say that $\hh\in \overline \F$ is a {\em generalized time derivative} of
$\hf\in \overline \F$, if for $T=T_0$,
\beq\label{gen derivative}
& &\lim_{\sigma\to 0+}
\left\|\frac{\hf(\cdotp+\sigma)-\hf(\cdotp)}\sigma-\hh
  \right\|_{\overline \F}=0,
\eeq
and write $\hh={\mathbb D} \hf$, or simply $\hh=\p_t\hf$.
We also need spaces with $s$ generalized derivatives,
$\F^s={\mathcal D}({\mathbb D}^s)$, with $s  \in \Bbb{Z}_+$,
and $\F^{\infty} = \bigcap_{s\in \Bbb{Z}_+}  \F^s$.
As in Remark \ref{remark 1a},
if $\hf\in \F^s$,
\beq\label{smoothness of gen. wave}
\omega^\hf_t\in \bigcap_{j=0}^s \big(C^{s-j}([T_0,\infty[,
{\mathcal D}(\M_e^j))
\cap \hbox{Ran}\,(\M_e)\big),
\eeq
so that
$\omega^\hf_t(T)\in {\bf H}^s_{\rm loc}(M^{\rm int})$ for $T\geq T_0$.

We need also the dual of the space  ${\mathcal D}(\M_e^s)$. Since
${\bf H}_0^s\subset {\mathcal D}({\mathcal M}_{\rm e}^s)$,
we have $({\mathcal D}({\mathcal M}_{\rm e}^s))'\subset
{\bf H}^{-s}$. Similarly,
${\bf H}^{-s}_0\subset ({\mathcal D}({\mathcal M}_{\rm e}^s))'$.
These facts will be used later to construct focusing sequences.

\subsection{Continuation of the boundary data}

Theorems \ref{blacho} and \ref{global control th}
make possible  to continue the boundary data, originally given
for $t \leq T$ to larger times $t>T$,
when $T$ is large enough, by using essentially the same ideas as in the
scalar case,
\cite{KKL,KL2} (see also \cite{BI} for another continuation method).

\begin{lemma}
\label{continuation} Admittance map
${\mathcal Z}^T$, given for  $T > 2\,{\rm rad}(M)$,
uniquely determines ${\mathcal Z}^t$ for any $t >0$.
\end{lemma}

{\em Proof:}
Let $2\e = T-2\,{\rm rad}(M)$. For
$f \in C^{\infty}_0([0, T], \,\Omega^1 \p M)$, Theorem
\ref{global control th} guarantees that
there is a sequence
$f_n \in C^{\infty}_0([\e, T], \,\Omega^1 \p M)$ with
\beq
\label{approx}
\lim_{n\to \infty} \omega^{f_n}_t(T)= \omega^f_t(T) \quad \hbox{in} 
\,\, L^2(\Omega^1M) \times
L^2(\Omega^2M),
\eeq
or, equivalently, in terms of the energy of the field,
\beq\label{Aa}
\lim_{n\to \infty}{\mathcal E}^T(\omega^{g_n})= 0, \quad  g_n = \p_t(f-f_n).
\eeq
{\newtext By using Theorem \ref{blacho}
one can verify  for an arbitrary sequence $(f_n)_{n=1}^\infty$
whether the
convergence condition (\ref{Aa}) is valid or not. Moreover,
the condition (\ref{Aa}) is valid for some sequence $(f_n)_{n=1}^\infty$.
Thus, when map ${\mathcal Z}^T$ is given, one can find
a sequence $(f_n)_{n=1}^\infty$ that satisfies condition (\ref{Aa}).
}

 From the definition (\ref{weak solution}) of a weak solution,
(\ref{approx}) implies
that
\beq
\label{conv}
\lim_{n\to \infty} {\bf n}\p_t(\omega^{f_n})^2|_{ \p M \times ]T,
\infty[}={\bf n}\p_t(\omega^f)^2|_{ \p M \times ]T, \infty[}.
\eeq
{\newtext Let $h_n(x,t)=f_n(x,t+\e)\in C^{\infty}_0([0, T-\e])$.
Since the function $\cZ^T h_n$
determines
  ${\bf n}\p_t(\omega^{f_n})^2
|_{ \p M \times ]T, T+\e[}$, we see that
$\cZ^T$  determines}
${\bf n}(\omega^f)^2|_{ \p M \times ]T, T+\e[}$.
Iterating this procedure, we construct
$\cZ^t$ for any $t >0$.
\hfill$\Box$

In sequel, we need $\cZ^T$ with various values
$T>2\hbox{rad}(M)$.
Taking into account Lemma \ref{continuation}, we denote simply
by $\cZ$ the admittance map
known for all $t$.

\begin{remark}
{\rm The controllability results, Theorem \ref{local control th}
together with the Blachovestchenskii
formula, Theorem \ref{blacho} make it possible to verify from
the knowledge of ${\mathcal Z}^T$ whether the condition $T>2\,{\rm rad}
(M)$ holds or not. Indeed, $T\leq 2\,{\rm rad}(M)$ if and only if
for any $\varepsilon>0$,
\[
  M(\partial M, (T-\epsilon)/2)\neq M(\partial M,T/2).
\]
{\newtext This is equivalent to that
there are $f_n\in C^\infty_0(\p M\times ]0,T/2[)$, $n=1,2,\dots$
such that $(\omega^{f_n})^1(T/2)$ form a Cauchy sequence
in $L^2(M)$, have norm one and converge to a function
that is orthogonal to all $(\omega^h)^1((T-\epsilon)/2)$,
  $h\in C^\infty_0(\p M\times ]0,(T-\epsilon)/2[)$.
When  ${\mathcal Z}^T$ is given,
this can be verified for all $f_n$ and $h$.}}
\end{remark}

\section{Inverse problem}

This chapter is devoted to the inverse problem of electrodynamics. Building
on properties of Maxwell's equations obtained in Chapter 1, we prove
the following uniqueness result.

\begin{theorem}\label{ip}
The boundary  $\partial M$ and the admittance map
${\mathcal Z}^T$, $T>2\,{\rm rad}(M)$,
uniquely determine  the  manifold
$M$, the travel metric $g$, and the scalar wave impedance $\alpha$.
\end{theorem}

The proof of this theorem consists of
several steps. The first
  is to reconstruct the Riemannian manifold $(M,g)$. Having $(M,g)$,
we then identify boundary sources which generate electromagnetic waves
focusing in a fixed point in $M$  at time $T> 2\,{\rm rad}(M)$.
These sources are instrumental
in reconstructing the impedance $\alpha$.
{\newtexx What is more, in section 
\ref{subs: data part} we prove a generalization of main Theorem
\ref{ip} for the case where the admittance map is given only a 
part of the boundary.}

 At the end, we return to $\R^3$ to
characterize group of transformations  of the parameters $\epsilon$
and $\mu$ leaving the boundary data intact.

\subsection{Reconstruction of the manifold}
\label{sect: reconstruction of manifold}

In this section we
determine the manifold $M$ and the travel
time metric $g$ from  the admittance map ${\mathcal Z}$.
The  idea is to use a slicing procedure to
control the supports of  the waves from the
boundary  in order to determine the
set of {\em boundary distance functions}.

We start by fixing certain notations.
Let   $T_0<T_1<T_2$ satisfy
\[
  T_0>2\,{\rm rad}(M),\quad T_1\geq T_0+{\rm diam}(M),
\quad T_2\geq 2 \,T_1.
\]
Let $\Gamma_j\subset\partial M$ be arbitrary open
disjoint sets and
$0<\tau_j^-<\tau_j^+<{\rm diam}(M)$ be arbitrary times, $1\leq j\leq J$.
We define a set $S = S(\{\Gamma_j,\tau_j^-,
\tau_j^+\}_{j=1}^J)\subset M$  as an intersection
of slices,
\begin{equation}\label{def of S}
   S =\bigcap_{j=1}^J  \left( M(\Gamma_j,\tau_j^+)\setminus
M(\Gamma_j,\tau_j^-) \right).
\end{equation}
(See Figure 2.)
Our first goal is to find, in terms of ${\mathcal Z}$,
whether the set $S$ contains an open ball or not.
To this end, we use the following definition.

\begin{definition}\label{support sources}
The set $Q= Q(\{\Gamma_j,\tau_j^-,
\tau_j^+\}_{j=1}^J) \subset {\mathcal F}^\infty $ consists
of  generalized sources
$\hat f$ such that the
waves, $\omega_t=\omega^{\hf}_t$,
satisfy
\begin{enumerate}
\item[(i)] $\omega_t^1(T_1)\in X(\Gamma_j,\tau_j^+)$, for
all $j$, $1\leq j\leq J$,
\item[(ii)] $\omega_t^2(T_1) = 0$,
\item[(iii)] $\omega_{tt}(T_1)=0$ in $M(\Gamma_j,\tau_j^-)$,
for
all $j$, $1\leq j\leq J$.
\end{enumerate}
\end{definition}

Observe that Maxwell's equations for $\omega_t=\omega^{\hf}_t$, $\hf\in Q$
imply that $\omega_{tt}=(0,\delta\omega_t^2,-d\omega_t^1,0)$, so in
particular at $t=T_1$, we have $\omega_{tt}(T_1)=(0,0,d\eta^1,0)$, where
$\eta^1 = -\omega_t^1(T_1)$ has the support property ${\rm supp}(d\eta^1
)\subset S$.

\begin{figure}[htbp]
\begin{center}
\psfrag{1}{\hspace{-1mm}$S$}
\psfrag{2}{\hspace{-2mm}$\Gamma_1$}
\psfrag{3}{$\Gamma_2$}
\psfrag{4}{\hspace{-3mm}$\tau_1^-$}
\psfrag{5}{$\tau_1^+$}
\includegraphics[width=5cm]{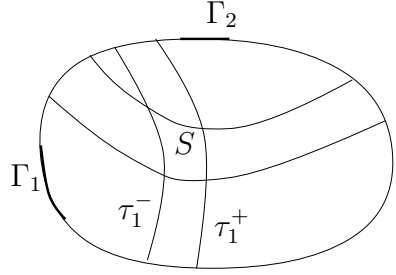} \label{pic 5}
\end{center}
\caption
{The set $S$ in the case when $J=2$.}
\end{figure}

The central tool for reconstruction of the manifold
is the following theorem.

\begin{theorem}\label{alternative}
Let $S$ and $Q$ be defined as above. Then:
\begin{enumerate}
\item If $S$ contains an open ball, then
${\rm dim}(Q)=\infty$,
\item If $S$ does not contain an open ball, then
$Q = \{0\}$.
\end{enumerate}
\end{theorem}

{\newtext The proof of Theorem \ref{alternative} is given later.}

\begin{theorem}\label{is in Z th}
For any $\hat f \in{\mathcal F}^\infty $ it is possible,
given ${\mathcal Z}$, to determine
whether $\hat f \in Q$ or not.
\end{theorem}

{\em Proof:} Let
$\hat f=(f_k)_{k=0}^\infty$ be a generalized source.
By  Remark \ref{Remark 1.},  Condition {\em (i)} of the definition of  $Q$
is equivalent to the existence of a sequence
$\hat h=(h_\ell)_{\ell=0}^\infty$,
$h_\ell \in\Cnull([0,\tau_j^+],\Omega^1\Gamma_j)$,
such that
\beq
\label{20.11.1}
  \lim_{k,\ell \to \infty} \| (\omega_t^{f_k})^1(T_1)
-(\omega_t^{h_\ell})^1(\tau_j^+)\| =0.
\eeq
By the linearity of the initial boundary value problem, we have
\[
\| (\omega_t^{f_k})^1(T_1)
-(\omega_t^{h_\ell})^1(\tau_j^+)\| =\|(\omega^{g_{k,\ell}})^1(T_1)\|,
\]
where the source $g_{k,\ell}$ is
\[
  g_{k,\ell}(t) = \partial _t\big(f_k(t) - h_\ell(t+\tau_j^+ - T_1)\big)
\in C_0^\infty
( ]0,T_1[,\Omega^1 \p M).
\]
Using Lemma \ref{blacho}, we can evaluate the norm
$\|(\omega^{g_{k,\ell}})^1(T_1)\|$ for various $(h_\ell)$
and thus verify the condition (\ref{20.11.1}).

In a similar fashion, Condition {\em (ii)}
is valid for $\hat f$, if
\[
  \lim_{k\to \infty}\|(\omega_t^{f_k})^2(T_1))\|= 0,
\]
which  can also be verified via ${\mathcal Z}$ by Lemma \ref{blacho}.

Finally, consider Condition {\em (iii)} for
$\hat f$ satisfying Conditions {\em (i)} and {\em (ii)}.
Observe that
\[
  (\partial_t+{\mathcal M})\omega_{tt}=0 \mbox{ in $M\times\R_+$},
\]
where $\omega_{tt}=\omega_{tt}^{\hat f}$,
and
$
  \bt\omega_{tt} =0 \mbox{ in $\partial M\times [T_0,\infty[$}.
$
If Condition {\em (iii)} holds, then, by the finite propagation speed,
\[
  \omega_{tt}=0\mbox{ in $K_j=\{(x,t)\in M\times\R_+\mid \tau(x,\Gamma_j)
+|t-T_1|<\tau_j^-\}$},
\]
i.e., $\omega_{tt}$ vanishes in the double cone of influence
of $\Gamma_j \times ]T_1- \tau_j^-, T_1+ \tau_j^-[$,
for all $j=1,\ldots,J$. Therefore, in each $K_j$,
$\omega_t$ does not depend on time, and, by  Condition {\em (ii)},
$\omega_t^2=0$ in $K_j$. Hence,
\begin{equation}\label{cond 3}
  \bn\omega_t^2= {\mathcal Z}f=0\mbox{ on $\Gamma_j\times
  ]T_1-\tau_j^-,T_1+\tau_j^-[$}.
\end{equation}
Conversely, if condition (\ref{cond 3}) holds together
with Conditions {\em (i)} and {\em (ii)}, then $\omega_t$ satisfies
\[
  (\partial_t+{\mathcal M})\omega_t =0\mbox{ in $M\times\R_+$}
\]
with the boundary conditions
\[
  \bt\omega_t^1=0,\quad \bn\omega_t^2=0\mbox{ in
  $\Gamma_j\times]T_1-\tau_j^-,T_1+\tau_j^-[$},
\]
because  $T_1-\tau_j^->T_0$,
so that $\hat{f}=0$ in $\Gamma_j\times]T_1-\tau_j^-,T_1+\tau_j^-[$.
Thus, by Theorem \ref{UCP}, $\omega_{tt}= 0$ in $K_j$ and,
in particular, Condition {\em (iii)} is
valid. As (\ref{cond 3}) is given in terms of
${\mathcal Z}$, this completes the proof.
\hfill$\Box$

{\em Proof of Theorem \ref{alternative}:} Assume that there is an
open ball $B \subset S$ and let
$\varphi \in \Omega^2M$ with $\hbox{supp}(\varphi) \subset \subset B$.
By Theorem \ref{global control th}, there is $\hat{f}  \in
\overline{\F}$ such that
\beq
\label{3a2}
\omega^{\hat{f}}_t(T_1) =(0, \delta \varphi,0,0),
\eeq
Clearly,  $\varphi \in {\mathcal D}^{\infty}({\mathcal M}_{ e})$, so that
$\hat{f} \in \F^{\infty} $.

Let us show that $\hat{f} \in Q$. Conditions {\em (i)}--{\em (ii)}
are immediate from (\ref{3a2}) and Theorem
\ref{local control th}. Finally, since
\[
\omega_{tt}^{\hat{f}}(T_1) = - {\mathcal M}\omega^{\hat{f}}(T_1)
=(0,0,d\delta \varphi,0),
\quad \hbox{supp}(d\delta \varphi) \subset\subset B,
\]
Condition {\em (iii)} is also valid. This proves the first part of the theorem.

To prove the second part, assume that $S$ does not contain
an open ball and, however, there is
$\hat{f} \in Q,\, \hat{f} \neq 0$. Let $\omega(t) =\omega^{\hat{f}}(t)$.
Then, by Conditions {\em (i)}--{\em (ii)},
\beq
\label{3a1}
\hbox{supp}(\omega_t(T_1)) \subset \bigcap_{j=1}^J M(\Gamma_j, \tau^+_j) = S^+,
\eeq
implying, due to $\omega_{tt}(T_1)=-{\mathcal M}\omega_t(T_1)$,
that ${\rm supp}(\omega_{tt}(T_1))\subset S^+$.
On the other hand, by Condition {\em (iii)},
\[
  \omega_{tt}(T_1)=0\quad\mbox{in $\bigcup_{j=1}^J
M(\Gamma_j,\tau_j^-)=S^-$.}
\]
Thus, ${\rm supp}(\omega_{tt}(T_1)) \subset S^+\setminus S^-$, which is
nowhere dense in $M$.
Since $\omega_{tt}(T_1)$ is smooth,
$\omega_{tt}(T_1)=0$, and, therefore,
\beq\label{refe +}
  d\omega_t^1(T_1)= -\omega_{tt}^2(T_1) =0.
\eeq
However, by Theorem \ref{global control th},
\[
  \omega_t^1(T_1) = \delta\eta^2 , \quad\mbox{with $\eta^2 \in
H(\delta, \Omega^2M)$}.
\]
Combining this equation with (\ref{refe +}) and using
$\bt\omega_t^1(T_1)=0$, we obtain, by Stokes' formula (\ref{stokes1}), that
\[
  (\omega_t^1(T_1),\omega_t^1(T_1))_{L^2} =
(\delta\eta^2,\omega_t^1(T_1))_{L^2}
  =(\eta^2,d\omega^1_t(T_1))_{L^2} = 0,
\]
i.e., $\omega^1_t(T_1)=0$.
Also, by Condition {\em (ii)}, $\omega^2_t(T_1) =0$. These
  imply that  $\hat f = 0$.

\hfill$\Box$

We are now ready to construct the set of  {\em boundary
distance functions}, $r_x$, which are defined, for any $x\in M$, as
continuous functions on $\p M$,
\[
  r_x: \p M\to\R_+,\quad r_x(z)=\tau(x,z), \quad z \in \partial M.
\]
They define {\it the
  boundary distance map} ${\mathcal R}:M\to C(\p M)$,
${\mathcal R}(x)=r_x$, which is continuous and injective,
\cite {Ku5,KKL}.
The set of all boundary distance functions, i.e., the image
of ${\mathcal R}$,
\[
  {\mathcal R}(M)=\{r_x\in C(\partial M)\mid x\in M\},
\]
can be endowed, in a natural way, with
a  differentiable structure and a metric tensor $\tilde g$,
so that $({\mathcal R}(M),\tilde g)$ becomes  isometric
to  $(M,g)$, see e.g. {\cite{Ku5,KKL}.
Hence, in order to reconstruct
$M$ (or more precisely, the isometry type of $M$),
it suffices to determine the set  ${\mathcal R}(M)$.
The following result is therefore crucial.

\begin{theorem}\label{thm: r(M)} 
For any $h\in C(\partial M)$,  it is possible, given
${\mathcal Z}$, to determine
whether $h\in {\mathcal R}(M)$ or not.
\end{theorem}

{\em Proof:} The proof is based on
a discrete approximation process.
First, we observe that
$h\in {\mathcal R}(M)$ if and only if, for any
finite subset
$\{z_1,\ldots,z_J\}$ of
$\p M$, there is an $x\in M$ with
\[
  h(z_j)= \tau(x,z_j),\quad 1\leq j\leq J.
\]
Denote $\tau_j = h(z_j)$. By the
continuity of the distance
function, $\tau:M\times \pM\to \R_+$,
the above condition is
equivalent to the following one:
For any $\varepsilon>0$, there are open sets $\Gamma_j\subset
\partial M$, $z_j \in \Gamma_j$  with ${\rm diam}(\Gamma_j)<\varepsilon$,
such that
\begin{equation}\label{int condition}
{\rm int}\bigg(\bigcap_{j=1}^J M(\Gamma_j,\tau_j+\varepsilon)
\setminus M(\Gamma_j,\tau_j-\varepsilon)\bigg)\neq
\emptyset.
\end{equation}
On the other hand, by Theorem \ref{alternative},
condition
(\ref{int condition}) is equivalent to
\[
  {\rm dim}\big( Q(\{\Gamma_j,\tau_j+\varepsilon,
\tau_j-\varepsilon\}_{j=1}^J)\big)=\infty,
\]
a condition that can be verified in terms of ${\mathcal Z}$
by means of  Theorem \ref{is in Z th}.
\hfill$\Box$

As a consequence, we obtain the main result of this section.

\begin{corollary}
The knowledge of the admittance map ${\mathcal Z}$
is sufficient for the reconstruction of the manifold $M$
and the travel time metric $g$.
\end{corollary}

Having found $(M,g)$, we proceed
to the reconstruction of the  impedance $\alpha$.

\subsection{Focusing sequences} \label{foc. sec}

In this section, we construct sequences
of  sources, $(\hat{f}_p)_{p=1}^\infty$ with the property that
$(\omega_{t}^{\hat{f}_p})^2(T_1)=0$ and
$\hbox{supp}\left((\omega_{t}^{\hat{f}_p})^1(T_1)\right)$
  converging,  when $p \to \infty$, to a single
point in $M^{\rm int}$, i.e., the time derivative of the electric
field focuses to a single point.

Let $y \in M^{\rm int}$ and $\ud_y$ denote the Dirac delta
at $y$ in the sense that
\[
  \int_M\frac 1\alpha\, \ud_y(x)\wedge *\phi(x)=\phi(y),
  \quad \hbox{for any $\phi\in C^\infty_0(M)$}.
\]
Since the Riemannian manifold $(M,g)$ is already found,
we can choose $\Gamma_{jp}\subset\partial M$,
$\,0<\tau_{jp}^-<\tau_{jp}^+<{\rm diam}(M)$, $j=1,\dots,J(p)$, so that
\beq
\label{20.11.3}
  S_{p+1}\subset S_p,\quad \bigcap_{p=1}^\infty S_p
  = \{y\},\quad S_p=S\big(\{\Gamma_{jp},\tau_{jp}^-,
\tau_{jp}^+\}_{j=1}^{J(p)}\big).
\eeq
Then,  $Q_p=Q\big(\{\Gamma_{jp},\tau_{jp}^-,
\tau_{jp}^+\}_{j=1}^{J(p)}\}\big)$ are spaces of boundary sources,
which correspond, by Definition \ref{support sources}, to the sets $S_p$.

\begin{definition}
For $y\in M^{\rm int}$, let $S_p, \, p=1,2,\dots,$ be given by
(\ref{20.11.3}). A sequence $(\hat f_p)_{p=1}^\infty$ with
$\hat f_p\in Q_p$, is called a {\em focusing sequence}
of boundary sources of order  $s$, $s \in \Bbb{Z}_+$, if
there is a  distribution form $A_y$ on $M$, $A_y \neq 0$, such that
\beq\label {fp-lim}
  \lim_{p\to \infty} (\omega^{\hat f_p}_t(T_1)
,\eta)_{{\bf L}^2}  = (A_y,\eta)_{{\bf L}^2}, \quad \hbox{when
$\eta\in {\mathcal D}({\mathcal M}_{\rm e}^s)$}.
\eeq
\end{definition}

With a slight abuse of notations, we use the same notation for the
inner product in ${\bf L}^2$ and for the distribution duality.
We denote a focusing sequences, converging to $y$,
by $\tilde f_y=(\hat f_p)_{p=1}^\infty$.

The following theorem characterizes a  class of limit
distributions
that can be produced by focusing sequences.  This class is large enough
for our further goal to solving inverse problem. What is
more, the sequences from this class can be {\newtext constructed} via 
the admittance map.

\begin{theorem}\label{1-forms}
(1)\  Let $y\in M^{\rm int}$ and $(\hat f_p)_{p=1}^\infty$ be a
sequence of boundary sources,
$\hat f_p  \in Q_p$. Given the admittance map,  ${\mathcal Z}$,
it is possible to determine,
for any $s \in \Bbb{Z}_+$, whether $(\hat f_p)$ is a focusing
sequence of order $s$ or not.

\noindent (2)\  Let  $\tilde f_y$ be a focusing sequence. Then ${\rm supp}(A_y)=\{y\}$.

\noindent (3)\  For $s=3$, the limit distribution $A_y$ has the form
\beq
\label{27.5}
A_y =(0,\delta(\lambda(y) {\underline{\delta}}_y),0,0),
\eeq
where $\lambda(y) \in \Lambda^2T^*_yM$.

\noindent (4)\  For any $y\in M^{\rm int}$ and $\lambda(y) \in \Lambda^2T^*_yM$,
  there is a focusing sequence $\tilde f_y$, of order $s=3$,
with $(A_y)^1= \delta(\lambda(y) {\underline{\delta}}_y)$.

\end{theorem}

{\em Proof:} {\it 1.} Take $\eta\in{\mathcal D}({\mathcal M} _{\rm e}^s)$
and decompose  it  as
$
  \eta=\eta_1+\eta_2,
$
where
\[
  \eta_1\in {\mathcal D}({\mathcal M}
_{\rm e}^s)\cap Y,\quad \eta_2\in
{\mathcal D}({\mathcal M}
_{\rm e}^s)\cap Y^\perp.
\]
As ${\mathcal U}(t)$ in (\ref{weak solution}) is unitary in
${\mathcal D}({\mathcal M}
_{\rm e}^s)$,
by  Theorem
\ref{global control th} there is a boundary source
$\hat h\in{\mathcal F}^s $ such that $\eta_1=
\omega_t^{\hat h}(T_1)$. Observe that $(\omega^{\hat f_p}_t(T_1),
\eta_2)_{{\bf L}^2} =0$,
so that
$(\hat f_p)$  is a focusing sequence if and only if there is a limit,
\begin{equation}\label{limit}
  (A_y,\eta)_{{\bf L}^2}=\lim_{p\to \infty} (\omega^{
\hat f_p}_t(T_1),\omega_t^{\hat h}(T_1))_{{\bf L}^2},  \quad \hbox{when
$\hat h\in{\mathcal F}^s $}.
\end{equation}
By Theorem \ref{blacho},  the existence of this
limit can be
verified in terms of ${\mathcal Z}$.

Conversely, assume that the limit (\ref{limit}) does
exist for all $\hat h\in{\mathcal F}^s$.
Then, by the Principle of Uniform Boundedness,
the functionals
\[
  \eta\mapsto(\omega^{\hat f_p}_t(T_1)
,\eta)_{{\bf L}^2}, \quad p \in \Bbb{Z}_+,
\]
are uniformly bounded  in the dual of
$\left({\mathcal D}({\mathcal M}_{\rm e}^s)\right)'$. By the
Banach-Alaoglu theorem, there is a weak$^*$-convergent
subsequence,
\[
\omega_t^{\hat f_p}(T_1)\to
  A_y\in \bigg({\mathcal D}({\mathcal M}_{\rm e}^s)
\bigg)',
\]
where $A_y$ is the sought distribution {\newtext for which
(\ref  {fp-lim}) is valid}.

{\it 2.} Let $\tilde{f}_y =(\hat f_p)_{p=1}^\infty$ be a focusing
sequence. Since $\hat f_p\in Q_p$, Condition {\em (ii)} of
Definition \ref{support sources} implies  that $A_y =(0,A_y^1,0,0)$ and
Conditions {\em (i)}--{\em (iii)}, together with (\ref{20.11.3}), yield
\beq
\label{28.1}
\hbox{supp} (dA_y^1)\subset
\liminf_{p\to \infty}\,
\hbox{supp}\big(d(\omega_t^{\hat f_p})^1 (T_1)\big) \subset
\bigcap_{p=1}^\infty S_p= \{y\}.
\eeq
As $\omega_t^{\hat f_p}(T_1) \in Y$,
$\,\delta A_y^1=\lim_{p\to \infty} \delta(\omega_t^{\hat f_p})^1 (T_1)=0$.
Thus,
\beq
\label{28.2}
\hbox{supp}(\Delta_{\alpha} A_y^1) \subset \{y\}.
\eeq
On the other hand,  by   Condition {\em (i)} of  Definition
\ref{support sources},
\ba
(\omega_t^{\hat f_p})^1(T_1) =0 \quad \hbox{in} \quad M \setminus S^+_p,\quad
  S^+_p =\bigcap_{j=1}^{J(p)}M(\Gamma_{jp},\tau_{jp}^+).
\ea
By the definition (\ref{fp-lim}) of a focusing sequence,
$A_y^1 = 0$ in $M \setminus S^+_p$.
  As $\hbox{rad}(M) < \hbox{diam}(M), $
we can always choose $\Gamma_{jp},\tau_{jp}^+$, so that $M \setminus S^+_p$
is non-empty.
By the unique continuation principle for elliptic systems (see e.g.
\cite{Isak1}), it then follows from the support property (\ref{28.2})
that
$
\hbox{supp}(A_y^1) \subset \{y\}.
$
Since $A_y$ is non-zero by assumption, $\hbox{supp}(A_y^1)= \{y\}$.

{\it 3.} Let $s=3$. {\newtext By part {\it 2.} of the theorem,
in local coordinates the components of
$A_y$ are finite sums of derivatives of delta-distributions.
Since}
$
A_y \in \left({\it D}({\mathcal M}^s_e) \right)' \subset {\bf H}^{-3}(M),
$
it follows that
\beq
\label{28.3}
A_y^1
= \sum_{i,j=1}^3 c_i^{j} \p_j{\underline{\delta}}_y \,dx^i +
\sum_{i=1}^3 \tilde{c}_i{\underline{\delta}}_y \,dx^i.
\eeq
Substituting (\ref{28.3}) into the
identity $\delta A_y^1=0$, we obtain (\ref{27.5}).

{\it 4.}
Let $\psi_p\in C^\infty_0(S_p)$, $p=1,2,\ldots$, be 2--forms that converge
to $\lambda {\underline{\delta}}_y$ in $H^{1-s}(\Omega^2M)$.
By the global control Theorem \ref{global control th},
there are boundary sources
$\hat f_p$ such that $\omega_t^{\hat f_p}(T_1)=(0,\delta \psi_p,0,0)$.
Then $\tilde{f}_y =(\hat f_p)_{p=1}^\infty$
is a desired focusing sequence.
  \hfill$\Box$
\smallskip

As $y$ runs over $M^{\rm int}$, we get a parametrized family of focusing
sequences $\{\tilde{f}_y\}_{y\in M^{\rm int}}$ which defines the map
$y\mapsto\lambda(y)$. However, the admittance map does not provide
a direct access to the values $\lambda(y)$. Although this mapping
is unknown, we have the following result.

\begin{lemma} \label{lem: 9.3}
  Given the admittance map ${\mathcal Z}$,
it is possible to determine
whether the map $y\mapsto\lambda(y)$ is  a
2--form valued $C^\infty$--function
in $M^{\rm int}$.
\end{lemma}

{\bf Proof:}
Let $\{A_y\}_{y\in M^{\rm int}}$ be a family of distributions
of the form (\ref{27.5}) corresponding to a family
$\{\tilde{f}_y\}_{ y \in M^{{\rm int}}}$ of focusing sequences.
Assume that $y\mapsto\lambda(y)$ is smooth, i.e., $\lambda\in\Omega^2
M^{\rm int}$.
Then, for any generalized source $\hat{h} =(h_j) \in {\mathcal F}^{\infty}$,
we have
\begin{eqnarray*}
  (A_y,\omega^{h_j}(T_1))_{{\bf L}^2} &=&
(\delta\lambda\ud_{y},(\omega^{h_j})^1
  (T_1))_{L^2} =(\lambda\ud_{y},d(\omega^{h_j})^1
  (T_1))_{L^2}\\
\noalign{\vskip4pt}
  &=& -(\lambda\ud_{y},(\omega^{h_j}_t)^2
  (T_1))_{L^2}.
\end{eqnarray*}
By taking the limit $j\to \infty$ of both sides and using the notation
$\langle \lambda,\eta \rangle_y=*(\lambda \wedge *\eta)$ for inner product
of
$\lambda, \eta \in \Lambda^kT^*_{y}M$, we arrive at the identity
\begin{equation}\label{4a1}
  (A_{y},\omega^{\hat{h}} (T_1))_{{\bf L}^2} =
  -\langle \lambda(y), ( \omega_t^{\hat{h}})^2(y,T_1) \rangle_{y}.
\end{equation}
{\newtext As
$
(A_{y},\omega^{\hat{h}}
  (T_1))_{{\bf L}^2}=
\lim_{j \to \infty}
\lim_{p\to \infty}
( \omega^{\hat{f}_{p,y}}_t(T_1),\omega^{h_j}(T_1))
_{{\bf L}^2},
$
we can evaluate (\ref{4a1})
in terms of ${\mathcal Z}$ by Theorem \ref{blacho}.}

Conversely, if $(A_{y},\omega^{\hat{h}}
  (T_1))_{{\bf L}^2} \in C^{\infty}(M^{{\rm int}})$  for any $\hat{h}
\in {\mathcal F}^{\infty}$,
then $\lambda (y) \in \Omega^2(M^{{\rm int}})$. Indeed, by
Theorem \ref{global control th},
for any $\varphi \in \Omega^1M$, ${\rm supp}(\varphi) \subset
\subset M^{{\rm int}}$,
there is a generalized boundary source $\hat{h} \in{\mathcal F}^{\infty}$
with $\omega_t^{\hat{h}}(T_1) = (0,0, -d\varphi,0)$,
and, by (\ref{4a1}),
\beq
\label{4a2}
  (A_{y},\omega^{\hat{h}}
  (T_1))_{{\bf L}^2}= \langle \lambda(y), d \varphi(y) \rangle_y
\in C^{\infty}(M^{{\rm int}}).
\eeq
As $\varphi$ is arbitrary, condition (\ref{4a2}) is equivalent to
that $\lambda(y)$ is
$C^\infty$--smooth
in $M^{\rm int}$.
\hfill$\Box$

Returning to (\ref{4a1}), we conclude that a focusing
  sequence $\{\tilde{f}_y\}$
gives rise  to a functional of $(\omega^{h_j}_t)^2(T_1))$.
It depends only on the value
of $(\omega^{h_j}_t)^2(T_1))$ at the point $y$ and will be called
{\em point evaluation functional} in the sequel. By the above result this
functional is determined upto an unknown factor $\lambda(y)$.
Hence,
by using three proper focusing sequences, we can evaluate the
2--form $(\omega_t^{\hat h})^2$ at any point in $M^{\rm int}$, up
to a linear transformation. The possibility to control
the precise form of this transformation
is discussed in the next section.

\begin{lemma}\label{gauge2}
Let $t>T_1$ and $\hh\in \overline {\F}$.
Given the admittance map  ${\mathcal Z}$,
it is possible to find the 2--forms
\beq\label{Bb}
K(y)(\omega^\hh_t(y,t))^2, \quad  y\in M^{{\rm int}}.
\eeq
Here  $K(y): \Lambda^2 T_y^*M\to \Lambda^2T_y^*M$ is
a smooth section of
${\rm End}(\Lambda^2T^*M^{{\rm int}})$ having
  maximal rank.
\end{lemma}

{\em Proof:} Let $U$ be a relatively open coordinate patch in $M$ with
2--forms $\xi_k\in \Omega^2 U$, $k=1,2,3$, linearly independent at any
$y\in U$. If $\{\tilde f_k(y)\}_{y\in U}$, $k=1,2,3$ are three families
of focusing sequences with the corresponding limiting 2--forms
$\lambda_k(y)$, we define the endomorphism $K_U(y)$ by
\beq
\label{25.11.1}
K_U(y)\omega^2(y) =\sum_{k=1}^3 \langle \lambda_k(y), \,
\omega^2(y) \rangle_y
\, \xi_k(y),\quad y\in U.
\eeq
{\newtext As we can
evaluate inner products (\ref{4a1}) by using Theorem \ref{blacho},
it is possible,  for any given
three families of focusing sequences $\{\tilde f_k(y)\}_{y\in U}$, $k=1,2,3$
and $\hat h$, to construct $K(y)(\omega^\hh_t(y,t))^2$ for $y\in U$,
$t>T$.} Further considerations are based on the result that we formulate
separately for future reference.

\begin{proposition}\label{endo}
Let $U\subset M^{\rm int}$ be open and $\xi_k\in\Omega^2 U$, $k=1,2,3$,
linearly independent at each $y\in U$. There are focusing sequences
$\{\tilde f_k(y)\}_{y\in U}$ such that the corresponding endomorphism
(\ref{25.11.1}) is $K_U(y)= I_y$, $y\in U$, the identity in
$\Lambda^2T_y^*M$.
\end{proposition}

{\em Proof:} Let $\lambda_k(y)$, $k=1,2,3$ form the dual basis of
$\xi_k(y)$, $k=1,2,3$,
\[
  \langle\lambda_k(y),\xi_\ell(y)\rangle_y =\delta_{k\ell}.
\]
It is a consequence of Theorem \ref{1-forms} that there are focusing
sequences $\tilde f_k(y)$ giving rise to the 2--forms $\lambda_k(y)$,
which shows the claim.
\hfill$\Box$

{\em End of the proof of Lemma \ref{gauge2}:}
By the above proposition, there are, for given linearly independent
$\xi_k(y)$, focusing sequences $\tilde f_k(y)$ so that $K_U(y)$
is of maximal rank. Moreover, since $K_U\omega^2(y)$ can be evaluated
for any $\omega^2 = (\omega_t^{\hat h})^2(T_1)$, the maximality of
the rank of $K_U(y)$ can be verified via ${\mathcal Z}$.

Let $U_j$, $j=1,\ldots,J$, be a finite covering of $M$ by coordinate
patches and $K_j$ the corresponding local endomorphism of the form
(\ref{25.11.1}) in $U_j\cap M^{\rm int}$.
{\newtext As we can compute (\ref{25.11.1}) for all $\hat h
\in {\mathcal F}^\infty$, $t>T$ and $x\in U$,
it possible to
verify that $K_j(y)=
K_\ell(y)$ for $y\in U_j\cap U_\ell$ for all $j$ and $\ell$. As by Proposition
\ref{endo} there are families of focusing sequences
for which this is true, we can construct the desired endomorphism.}
\hfill$\Box$

\subsection{Reconstruction of the wave impedance}
\label{sect: wave impedance}

So far, we have found the waves $(\omega^{\hat h}_t)^2(t)$,
$t>T_1$, up to a linear transformation $K$ which, at this stage, is unknown.
Since the choice of the focusing sequences is non-unique, we will
choose them in such a manner that the endomorphism $K$ becomes as
simple as possible, i.e., $K = c_0I$, an identity up to a constant
multiplier. The first step in this direction is to consider the
polarization of the {\em electric Green's function}, defined as
the solution of the following initial boundary value problem,
\begin{eqnarray}\label{green equations}
  (\partial_t + {\mathcal M})G_{\rm e} (x,y,t)&=& 0\mbox{ in $M\times \R_+$},
\nonumber
\\
\noalign{\vskip4pt}
  \bt G_{\rm e}(x,y,t) &=& 0\mbox{ in $(x,t)\in \partial M\times \R_+$},
\\
\noalign{\vskip4pt}
G_{\rm e}(x,y,t)|_{t=0} &=& (0,\delta(\lambda \ud_y),0,0).\nonumber
\end{eqnarray}
{\newtexx Sometimes, we denote $G_{\rm e} (x,y,t)=G_{\rm e} (x,y,t;\lambda)$
to indicate the source $\lambda\in \Lambda^2 T^*yM$.}
Assume that $\tilde h=\tilde h_y$ is a focusing sequence that produces a
wave focusing at $y$, the corresponding 2--form being $\lambda$. Since
the boundary sources are off when $t>T_1$, we must have
\begin{equation}\label{el green}
G_{\rm e}(x,y,t) =
\omega_{t}^{\tilde h}(x,t + T_1).
\end{equation}
{\newtext On the other hand, by
Lemma \ref{gauge2}, we can calculate
the 2--forms $K(x)(\omega_{t}^{\tilde h})^2(x,t+T_1)$
for $x\in M^{\rm int}$ and $t>0$.}
Hence, we know the electric Green's function up to a linear transformation.

Let us denote by $\Phi =\Phi^\lambda(x,y,t)$
the standard Green's 2--form, satisfying
\begin{eqnarray}
  (\partial_t^2+\Delta^2_{\alpha}) \Phi(x,y,t) &=&0
\hbox{ in  $M\times \R_+$}, \nonumber \\
\noalign{\vskip4pt}
\Phi(x,y,t)|_{\pM \times \R_+} &=& 0,\\
\noalign{\vskip4pt}
\Phi(x,y,0)=0,& & \Phi_t(x,y,0)= \lambda(y)
\ud_y(x),
\label{1-form Green's function}
\end{eqnarray}
where $\lambda(y) \in  \Lambda^2T^*_yM$ and
the boundary condition in (\ref{1-form Green's function}) means that
all three components of  $\Phi$ vanish on $\pM \times \R_+$.

Let $\widetilde{G}_{{\rm e}}= \widetilde{G}_{{\rm e}}(x,y,t)$
be defined as
\beq
\label{6a1}
  \tilde{G}_{\rm e} =(\partial_t-{\mathcal M})
(0,\delta \Phi,0,0)
  = (0,\partial_t \delta \Phi, -d\delta\Phi, 0).
\eeq
As
$
  (\partial_t^2 + {\bf \Delta}_{\alpha})=
  (\partial_t +{\mathcal M})(\partial_t -{\mathcal M}),
$
$\tilde{G}_{\rm e}$ satisfies the complete Maxwell system and,
by (\ref{1-form Green's function}),  the initial condition in
(\ref{green equations}). By  the unit propagation speed,
$ \tilde{G}_{{\rm e}}  =0$ near $\p M \times ]0,\tau(y,\pM)[$,
satisfying  the boundary condition in (\ref{green equations}). Thus,
$\tilde{G}_{\rm e}(x,y,t) =G_{\rm e}(x,y,t)$ for $t<\tau(y,\pM)$.

To further the study of $G_{\rm e}$, we formulate the following
result proved in the Appendix.

\begin{lemma}\label{asymptotics of green}
For every $y\in M^{\rm int}$ there is an open neighborhood
$U\subset M^{\rm int}$ of $y$, a positive $t_y$ and a mapping
$Q_y(x)$ that is smooth with respect to $x\in U$, where
$Q_y(x):\Lambda^2T_y^* M\to\Lambda^2T_x^* M$ is bijective,
such that
\beq
\label{7a1}
\Phi(x,y,t)=Q_y(x)\lambda\ud(t^2-\tau^2(x,y))+r(x,y,t),\quad
(x,t)\in U\times ]0,t_y[.
\eeq
Moreover, with some smooth $Q_y^p(x):\Lambda^2T^*_yM\to
\Lambda^2T^*_xM$, $p=1,2$  and $C^{1,1}$--smooth 2--form $\hat r(x,y,t)$,
the remainder can be written as
\begin{equation}\label{remainder}
  r(x,y,t) = \sum_{p=1,2}Q_y^p(x)\lambda(t^2-\tau^2(x,y))_+^{p-1} +
\hat r(x,y,t).
\end{equation}
\end{lemma}

By (\ref{6a1}), it follows from (\ref{7a1}) that, for sufficiently
small $t$,
\beq
\label{7a2}
G_{{\rm e}} = (0,G_{\rm e}^1,G_{\rm e}^2,0) + r_1,
\eeq
where
\begin{eqnarray*}
G_{\rm e}^1 = -2t*(d\tau^2\wedge *Q_y\lambda) \ud^{(2)}
(t^2-\tau^2),\ \
G_{\rm e}^2 =d\tau^2 \wedge *(d\tau^2\wedge *Q_y\lambda)\ud^{(2)}
(t^2-\tau^2),
\end{eqnarray*}
and $r_1$ is linear combination of
a bounded function, the
delta distribution on $\partial B_y(t)$ and its first derivative,
$B_y(t)$ being the ball of radius $t$ centered in $y$.
Using  Lemma \ref{gauge2} and (\ref{el green}) together with (\ref{7a2}),
we obtain the following result.

\begin{lemma}\label{values of Maxwell}
Given the admittance map ${\mathcal Z}$, it is
possible to find the distribution 2--form
\ba
K(x)G_{\rm e}^2(x,y,t)= K(x)(\omega^{\tilde h}_t)^2(x,t+T_1),
\ea
where $K \in {\rm End}(\Omega^2 M^{\rm int})$ and $t>0$.
Moreover, the leading singularity of this form
when $0<t<t_y$ determine the 2--form
\beq\label{wave front value}
  K(x) \big(d\tau^2(x,y) \wedge *(d\tau^2(x,y)\wedge *Q_y(x)\lambda)\big),\quad
x\in\partial B_y(t).
\eeq
\end{lemma}
The linear transformation $K(x)$ of Lemma \ref{gauge2} depends on
$\tilde{f}_k(x)$, $\xi_k(x)$,  $k=1,2,3$. Our next goal is to formulate
conditions, in terms of  ${\mathcal Z}$, on
$\tilde{f}_k,\xi_k$ to make $K$ isotropic, i.e.,
\beq\label{K-requirement}
K(x)=c(x)I, \quad c \in C^{\infty}(M^{{\rm int}}), \quad c(x) \neq 0.
\eeq
To this end, observe that for $ \lambda \in \Lambda^2T^*_yM,$
\beq
\label{21.11.2}
{\bf t}_{B_y(t)} \big(d\tau^2 \wedge *(d\tau^2\wedge *Q_y\lambda)\big) =0,
\eeq
  where ${\bf t}_{B_y(t)} \omega^k$
is the tangential component of $\omega^k$ on $\p B_y(t)$.
Physically, condition (\ref{21.11.2}) corresponds  to the 
orthogonality of  the polarization
of the magnetic flux density and the direction of wave propagation.
(See Figure 3.)
If $K$ is isotropic, we have
\beq
\label{orthog}
{\bf t}_{B_y(t)} \bigg(K\big(d\tau^2 \wedge *(d\tau^2\wedge *
Q_y\lambda)\big)\bigg) =0.
\eeq

\begin{figure}[htbp]
\begin{center}
\psfrag{1}{$\vec v$}
\psfrag{2}{$\vec w$}
\includegraphics[width=10cm]{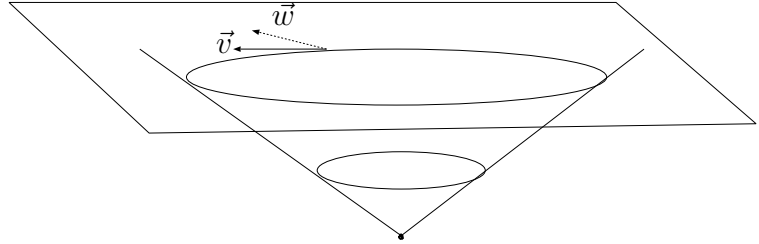} \label{pic 6}
\end{center}
\caption{Vector $\vec v$ is the right polarization of the electromagnetic wave
in the plane
$M\times \{t\}$. The reconstructed polarization $\vec w$ has
wrong direction, if the transformation matrix $K(x)$ is not isotropic.
}
\end{figure}

Conversely, we show that condition (\ref{orthog})
{\newtext for all $y\in M^{\rm int}$ and $t<t_y$} guarantees that
$K$ is isotropic. What is more,
the condition (\ref{orthog}) is verifiable from the knowledge of
${\mathcal Z}$.
Indeed, for  $\lambda(y) \in \Lambda^2T^*_yM$ and  $t=\tau(x,y)$,
(\ref{orthog}) means that
$K\big(d\tau^2 \wedge *(d\tau^2\wedge *
Q_y\lambda)\big)$ is normal to $T_x\partial B_y(t)\subset T_x M$, i.e,
for vectors $X$, $Y\in T_x\partial B_y(t)$, we have
\[
K(x)\big(d\tau^2(x,y) \wedge *(d\tau^2(x,y)\wedge *
Q_y(x)\lambda(y))\big)(X,Y)=0.
\]
Observe that when $\lambda(y)$ runs through $\Lambda^2 T_y^*M$,
then $*(Q_y(x)\lambda(y))$ runs through $T^*_x M$.
Now we may vary $y$ and $t$ with fixed $x$ such that $\tau(x,y)=t$,
making $T_x\partial B_y(t)$ run through the Grassmannian manifold
$G_{3,2}(T_xM)$. Transformation $K(x)$ is kept invariant in this
variation. Hence, we deduce that $K(x)$ keeps any 2--dimensional
subspace of $\Lambda^2T_x^*M$ invariant, so it must be isotropic
as claimed.

Assume that the
  the focusing sequences used for the point evaluation functionals
are chosen so that $K(x)=c(x)I$. For any generalized source $\hat f
\in {\overline{\mathcal F}}$, we may thus evaluate
\ba
(\tilde \omega_t^{\hat{f}})^2(x,T_1)=c(x) (\omega_t^{\hat{f}})^2(x,T_1),
\ea
with yet unknown $c(x)$.
Since $\omega^{\hat f}$ satisfies Maxwell's equations, we have
\[
  d(\tilde \omega_t^{\hat{f}})^2 = dc\wedge(\omega_t^{\hat{f}})^2
  + cd(\omega_t^{\hat{f}})^2 = dc\wedge(\omega_t^{\hat{f}})^2.
\]
The global control Theorem \ref{global control th} thus asserts that
$c(x)=c_0$ is equivalent to
\beq
\label{const}
d( \tilde \omega_t^{\hat{f}})^2(x,T_1)=0, \quad
\hbox{for all $\hat f\in {\overline{\mathcal F}}$},
\eeq
a condition that is  verifiable from the knowledge of
${\mathcal Z}$. Hence, the focusing sequences used for point evaluation
can be chosen such that $c(x)=c_0\neq 0$.

To proceed with reconstruction of $\alpha$, consider
the inner product,
\begin{equation}\label{inner product}
\int_M (\tilde \omega_t^{\hat{f}})^2(x,T_1)\wedge *
(\tilde \omega_t^{\hat{h}})^2(x,T_1)=
c_0^2 \int_M
  (\omega_t^{\hat{f}})^2(x,T_1)\wedge *
(\omega_t^{\hat{h}})^2(x,T_1),
\end{equation}
which can be found via ${\mathcal Z}$. On the other hand, by Theorem
\ref{blacho}, ${\mathcal Z}$ determines the energy inner product,
\ba
\frac 12 \int_M \frac 1{\alpha(x)}(\omega_t^{\hat{f}})^2(x,T_1)\wedge *
(\omega_t^{\hat{h}})^2(x,T_1).
\ea
By choosing the boundary source $\hat h = \hat h_j$ such that
$\tilde h =(\hat h_j)_{j=1}^\infty$ is a focusing sequence and by comparing the
above inner products at the limit $j\to \infty$, we recover the value
$c_0^2\alpha(x)$ at any point $x\in M$.

Finally, we notice e.g. by considering the
energy integrals that the admittance map has the
scaling property  ${\mathcal Z}_{(M,g,c_0^2\alpha)}=
c^{-2}_0 \,{\mathcal Z}_{(M,g,\alpha)}$, with evident notations.
Therefore, given  ${\mathcal Z}$ and $(g,c_0^2\alpha)$ already reconstructed,
it is also possible to determine $c_0$ and
hence $\alpha$.
This completes  the proof of Theorem \ref{ip}.
\hfill$\Box$

\subsection{Data given on the part of the boundary}
\label{subs: data part}

In this section, we generalize the proof of Theorem
\ref{ip} for the case when the data is given on a non-empty
open subset $\Gamma\subset \p M$.
In this case, instead of the complete admittance map
${\mathcal
Z}^T$ we are given the local admittance map ${\mathcal
Z}^T_\Gamma$, defined by
\[
{\mathcal  Z}^T_\Gamma f= 
{\mathcal Z}^T f|_{\Gamma\times ]0,T[},\quad f\in \Cnull([0,T],\Omega^1\Gamma),
\]
where $\Omega^1\Gamma$ is the space of 1-forms
$f\in \Omega^1\p M$ supported on $\Gamma$. Denote
${\mathcal Z}={\mathcal Z}^T$ with $T=\infty$
{\newtexxt and recall that ${\rm rad}_\Gamma(M)$ is the geodesic radius of $M$ with respect
to $\Gamma$, see (\ref{23.6.1}).
}

\begin{theorem}\label{ip2}
Given $\Gamma$, the local admittance map
${\mathcal Z}^T_\Gamma$, $T>2\, {\rm rad}_\Gamma(M)$,
uniquely determines  the  manifold
$M$, the  metric $g$, and the scalar wave impedance $\alpha$.
\end{theorem}

{\bf Proof.} Here we use notations of section \ref
{generalized sources section}.
By Theorem
\ref{global control th} we have that the set
 ${\mathcal
F}_\Gamma=C^\infty_0(]0,T_0[, \Omega^1\Gamma )/\hspace{-2mm}\sim$
with  $T>T_0>2
\,\hbox{rad}_\Gamma(M)$ is a dense subset of
 $\overline {\mathcal F}$. Thus we can identify
 $\overline {\mathcal F}_\Gamma$ with $\overline {\mathcal F}$.
This makes it possible to use,
 when the data is given on a part of the boundary,
 all the results about generalized sources
obtained in section \ref{generalized sources section} for the whole boundary.
In particular, we can define sets  ${\mathcal
F}^s_\Gamma\subset  \overline {\mathcal
F}_\Gamma$ that can be identified with  ${\mathcal
F}^s$.

Exactly as in section
\ref{sect: continuation} we can show that the
local admittance
map ${\mathcal Z}^T_\Gamma, \, T>2
\,\hbox{rad}_\Gamma(M),$
determine the map  ${\mathcal Z}^t_\Gamma$
for
all $ t>0$, i.e.,  the map  ${\mathcal Z}_\Gamma$.

Our first aim is to reconstruct $(M,g)$ near $\Gamma$.
For this, let $\theta_{{}_{\Gamma}}:\Gamma\to \R$ be
\ba
\theta_{{}_{\Gamma}}(z) =  \sup\{ s>0: \,
\tau(\gamma_{z,\nu}(s),\Gamma)=s \}
\ea
and
\ba
M_{\Gamma}=\{\gamma_{z,\nu}(s)\in M\ :\ z\in \Gamma,\
0\leq s<\theta_{{}_{\Gamma}}(z) \},
\ea
{\newwwtext where $\gamma_{z,\nu}(s)$ is the geodesic starting from $z \in \p M$ in the normal
direction.}

\begin{lemma}\label{part of bnd, alpha}
Given $\Gamma$, the local admittance map   ${\mathcal Z}_\Gamma$
determines the function $\theta_{{}_{\Gamma}}:\Gamma\to \R$, the Riemannian
manifold $(M_\Gamma,g)$, and the wave impedance $\alpha$ on $M_\Gamma$.
\end{lemma}

{\em Proof:}
Let  $z\in \Gamma$.
Using notations of section \ref{sect: reconstruction of manifold}
we see that $s\leq \theta_{{}_{\Gamma}}(z)$
if and only if $S=M(\Gamma_z,s)\setminus M(\Gamma,s-\e)$
has non-empty interior for all open 
$\Gamma_z\subset \Gamma$ containing $z$
and $\e>0$. In turn,
$S=(M(\Gamma_z,s)\setminus M(\Gamma_z,0))
\cap (M(\Gamma,s)\setminus M(\Gamma,s-\e))$ is of the form
(\ref{def of S}). Thus, using Theorem
\ref{is in Z th}
{\newtexxt with ${\mathcal Z}_\Gamma$ instead of ${\mathcal Z}$,}
 we can find out whether $S$ has non-empty interior
or not. Thus we can find $\theta_{{}_{\Gamma}}(z)$.

At this stage, we can proceed in a similar manner to the proof of
Theorem \ref{thm: r(M)} to find functions
$R_{\Gamma}(M_{\Gamma})=\{r_x|_\Gamma\,:\ x\in M_\Gamma\}$.
 To do this, we just use the procedure presented
in the proof of Theorem \ref{thm: r(M)} but consider only
functions $h\in C(\Gamma)$ that have a unique global minimum,
say $z_0\in \Gamma$, with $h(z_0)< \theta_{{}_{\Gamma}}(z_0)$.
After construction of this set, we see as in \cite[Sect. 4.4]{KKL} that the set
$R_{\Gamma}(M_{\Gamma})$ determines the Riemannian manifold
$(M_\Gamma,g)$.

Reconstruction of $\alpha$ in $M_\Gamma$ follows the same route
as with data given on the whole boundary by restricting
our attention to focusing sequences corresponding to points $y\in M_\Gamma$
and using identification of
${\mathcal F}^s_\Gamma$ with  ${\mathcal F}^s$.
\hfill$\Box$

In the next step we will show that we can find the admittance
map on the boundary of an arbitrary ball $B\subset M_\Gamma$.
We will denote by ${\mathcal Z}_{\p B}$ the local admittance
map defined by using the manifold $M\setminus B$ instead of $M$
and $\p B$ instead of $\Gamma$. For similar arguments in
the scalar case, see \cite{KKL2}.

\begin{proposition}\label{prop A9}
 Given $(M_\Gamma,g,\alpha)$ and
the map ${\mathcal Z}_\Gamma$ for $(M,g,\alpha)$, we can
find the local admittance map ${\mathcal Z}_{\p B}$
for $(M\setminus B,g,\alpha)$.
\end{proposition}

{\bf Proof.} First we observe that
 $(M_\Gamma,g,\alpha)$ and ${\mathcal Z}_\Gamma$
determine the values of the electric Green's function $ G^2_{\rm e}(x,y,t;\lambda)$ for any $x,y\in M_\Gamma$,
$t>0$, and $\lambda\in \Lambda^2T^*_yM$.
This result is proven in section  \ref{sect: wave impedance}
in the case when the admittance map is given
on the whole boundary and the proof can be directly extended
to the considered case.

Consider now the initial boundary value problem
\beq\label{ibvp eta}
& &\p_t \eta + {\mathcal M}\eta =\kappa, 
 \quad\hbox{in }M\times\R_+, \quad
\bt\eta\big|_{\partial M\times \R_+}=0,\quad 
\eta(0)=0, 
\eeq
where  $\kappa=(0,\delta \beta,0,0)$ with $\beta\in \Cnull(\R_+,
\Omega^2B)$ and denote its solution by $\eta=\eta_\beta=
(0,\eta^1,\eta^2,0)$.
Writing $\eta^2$ in terms of the electric Green's function
we obtain
\beq\label{eta 2}
\eta^2(x,t)=\int_{\R_+}\int_B
 G^2_{\rm e}(x,y,t-t';\frac{\beta(y,t')}{\alpha(y)})\, dV_g(y)dt',
\eeq
where $dV_g$ is the Riemannian volume measure 
on $(M,g)$.
Using equation (\ref{ibvp eta}) we also find 
\beq\label{eta 1}
\eta^1(x,t)=\int_0^t (\delta \eta^2(x,t')+\delta\beta(x,t'))dt'.
\eeq

{\newtexxt We continue  the proof with the following lemma;}

\begin{lemma}\label{lem D}
Let $\omega=(0,\omega^1,\omega^2,0)\in \Cnull(\R_+,\Omega M)$ be a solution
of Maxwell's equations
(\ref{MA1}) and (\ref{MF1}) in $(M\setminus B)\times
\R_+$ which satisfies the electric boundary condition
$\bt \omega=0$ on $\p M\times \R_+$ and
initial condition $\omega(0)=0$. Then there
is $\beta\in \Cnull(\R_+,
\Omega^2B)$ such that the solution    $\eta_\beta$
of initial boundary value problem (\ref {ibvp eta})
coincides with $\omega_{tt}$ in 
 $(M\setminus B)\times
\R_+$.\end{lemma}

{\bf Proof.} Let $\tilde \omega=(0,\tilde \omega^1,
\tilde \omega^2,0)\in \Cnull(\R_+,\Omega M)$
be an arbitrary  smooth continuation of $\omega$ into
$B\times \R_+$. Let
\ba
\rho=(0,\rho^1,\rho^2,0),\quad \rho^1=-\delta d\tilde\omega^1,
\quad \rho^2=- d \delta \tilde \omega^2.
\ea
Then $\rho=\omega_{tt} $ in $(M\setminus B)\times
\R_+$, and $\rho$ satisfies Maxwell's equations 
\ba
& &\rho_t^1-\delta\rho^2=\delta d a^1,\quad a^1=-\tilde \omega^1_t+\delta
\tilde \omega^2,\\
& &\rho_t^2+d\rho^1=d\delta a^2,\quad a^2=-\tilde \omega^2_t-d
\tilde \omega^1,
\ea
in $M\times \R_+$. Then $\eta=
(0,\rho^1-\delta a^2,\rho^2,0)$ satisfies 
 the initial boundary value problem (\ref{ibvp eta})
with $\beta=d a^1-a^2_t$ supported in $B\times \R_+$.
In particular, $\omega_{tt}=\eta_\beta$ in 
 $(M\setminus B)\times
\R_+$.
\hfill$\Box$

{\newtexxt To complete the proof of Proposition \ref{prop A9}
we start with an arbitrary $\beta\in \Cnull(\R_+,\Omega^2B)$ and find, using 
formulae 
(\ref{eta 2}), (\ref{eta 1}), the wave $\eta_{\beta}(x,t)$ for  $x \in M_{\Gamma}$.
 Let $\omega(x,t) $ be now defined as 
\beq
\label{23.6.2}
\omega^0(x,t)=0,\quad
\omega^1(x,t)=
\int_0^t\int_0^{t'} \eta_\beta^1(x,t'')\,dt''dt',\\
\nonumber
\omega^2(x,t)=
\int_0^t\int_0^{t'} \eta_\beta^2(x,t'')\,dt''dt', 
\quad
\omega^3(x,t)=0.
\eeq
Then $\omega(t) $ is a solution of  the initial-boundary value problem,
\ba
\omega_t + \mathcal{M} \omega =0 \quad \hbox{in} \,\,(M\setminus B)\times
\R_+,  \quad \omega(0)=0, \\
 \bt \omega|_{\p M \times \R_+}=0,
\quad
\bt \omega|_{\p B \times \R_+}=(0, f_{\beta}, -\int_0^t df_{\beta}),
\ea
where
\beq
\label{24.6.1}
f_{\beta} = \int_0^t\int_0^{t'}\bt \eta_\beta^1(x,t'')\,dt''dt'
\in  \Cnull(\R_+,\Omega^1\p B).
\eeq
Using again formulae (\ref{eta 2}), (\ref{eta 1}), we see that 
$(M_\Gamma,g,\alpha)$, 
${\mathcal Z}_\Gamma$ determine the map
\beq
\label{24.6.2}
f^1_{\beta} \longrightarrow \bn \omega^2|_{\p B \times \R_+} = 
\int_0^t\int_0^{t'}\bn \eta_\beta^2(x,t'')\,dt''dt'
\in \Cnull(\R_+,\Omega^1\p B)
\eeq
for any $f^1_{\beta}$ of form (\ref{24.6.1}).
As according to Lemma \ref{lem D}. the map $\beta \to f^1_{\beta}$ is a surjective map
from $\Cnull(\R_+,\Omega^2 B)$ onto $\Cnull(\R_+,\Omega^1\p B)$,
the map
(\ref{24.6.2}) determines $\mathcal{Z}_{\Gamma}$.
}
This proves Proposition \ref{prop A9}.
\hfill$\Box$

Having found ${\mathcal Z}_{\p B}$ we construct the 
Riemannian manifold $M_{\p B} \subset M\setminus B$,
metric $g$ on $M_{\p B}$ and impedance $\alpha$ on
$M_{\p B}$. Here $M_{\p B}$ is defined in a similar way as $M_\Gamma$
changing $M$ to $M\setminus B$ and $\Gamma$ to
$\p B$. Combining this with the previous results,
we  find the part $M_\Gamma\cup M_{\p B}$
of $M$ as well as the metric $g$ and the wave impedance $\alpha$
on it. Iterating this procedure, we reconstruct,
in finite number of steps, the whole manifold
$(M,g,\alpha)$. For detail, see \cite[Sect 4.4.9]{KKL}.
This proves Theorem \ref{ip2}.
\hfill $\Box$

\subsection{Inverse problem for Maxwell equations in $\R^3$}
 
In this section, the  uniqueness results
for Maxwell's equations on a  manifold are used to 
characterize the non-uniqueness of  
inverse problems for Maxwell's equations
(\ref{Maxwell-Faraday 1})--(\ref{Maxwell--Ampere 1}) in a 
bounded domain of $M \subset \R^3$
with the Euclidean metric $(g_0)_{ij}=\delta_{ij}$.

Let $M_j\subset\R^3$, $j=1,2$, be two 
bounded smooth closed domains
with a common part $\Gamma$ of their boundaries, $\Gamma\subset \p M_1\cap
\p M_2$. Let
$\epsilon_j$ and $\mu_j$, $j=1,2$ be permittivity
and permeability matrices in $M_j$, respectively, with
$\mu_j=\alpha^2_j\epsilon_j$, with $\alpha_j>0$ being the corresponding scalar
impedances. 
Assume that the local
admittance maps ${\mathcal Z}_{\Gamma,j}$ for $(M_j,\epsilon_j,\mu_j)$ coincide. 
{\newtext By Theorem \ref{ip2}, both  $(M_1,\epsilon_1,\mu_1)$
and  $(M_2,\epsilon_2,\mu_2)$ correspond
to the same abstract
manifold 
$(\tilde M,\tilde g,\tilde \alpha)$ which is  uniquely determined by 
${\mathcal Z}_{\Gamma,j}$ with the part $\Gamma$ corresponding to a part
 $\tilde \Gamma\subset \p \tilde M$. This implies that there are 
embeddings $F_j: \tilde M\to M_j\subset\R^3$ of the manifold $\tilde M$ in
the Euclidean space such that the metric tensors and
the wave impedances satisfy $\tilde g=(F_j)^*g_j$ and
 $\tilde \alpha =(F_j)^*\alpha_j$ and $F_1|_{ \tilde \Gamma}=F_2|_{ \tilde \Gamma}$.} 
{\newtexxt Recall that $g_j$ are determined by expression (\ref{metric}) with $\e_j$
and $ \mu_j$ in place of $\e, \mu$.}
Embeddings $F_j$ induce a diffeomorphism
\begin{equation}\label{19.0}
 \Phi = F_2\circ F_1^{-1}:M_1\to M_2,\quad \Phi\big|_{\Gamma}=
 {\rm id}.
\end{equation}
Consider two vector fields $X_1$ and $Y_1$ in $M_1$, and 
{\newtext denote
$X_2= D \Phi X_1$, $Y_2=D \Phi Y_1$.} The electric energy inner product
for the corresponding 1--forms $\omega^1$, $\eta^1\in\Omega^1 M$
is invariant, i.e., we have
\[
 \int_{M_1}g_0(X_1,\epsilon_1 Y_1)dV_0 =
\int_{\tilde M}\frac 1{\tilde \alpha}\omega^1
\wedge * \eta^1 =  \int_{M_2}g_0(X_2,\epsilon_2 Y_2)dV_0.
\]
On the other hand, as $X_2=D\Phi\, X_1$ and  $Y_2=D \Phi Y_1$, 
\[
\int_{M_2}g_0(X_2,\epsilon_2 Y_2)dV_0 = \int_{M_1}g_0(X_1,\Phi^*\epsilon_2
Y_1)dV_0,
\]
{\newtext where
\begin{equation}\label{19.3}
 \Phi^*\epsilon_2 = \frac 1{{\rm det}\,D\Phi} (D\Phi)^{\rm T}\, (\epsilon_2\circ\Phi)\,
D\Phi.
\end{equation}
Since $X_1$ and $Y_1$ are arbitrary, we must have 
$\epsilon_1=\Phi^*\epsilon_2$. Similar reasoning shows that
$\mu_1=\Phi^*\mu_2$. 
}

Thus we have proven the following result.

\begin{theorem} 
\label{group}
Let $M_1,M_2\subset \R^3$ be bounded smooth domains and
$\Gamma\subset \p M_1\cap \p M_2$ be  open and
non-empty. Let
 ${\mathcal Z}_{\Gamma,1}$ and ${\mathcal Z}_{\Gamma,2}$
be the local admittance maps corresponding
to $(M_1,\epsilon_1,\mu_1)$
and $(M_2,\epsilon_2,\mu_2)$, respectively. Then 
 ${\mathcal Z}_{\Gamma,1}={\mathcal Z}_{\Gamma,2}$
if and only if there is a diffeomorphism $\Phi:M_1\to M_2$,
$\Phi\big|_{\Gamma}={\rm id}$ and 
$\epsilon_1=\Phi^*\epsilon_2$, $\mu_1=\Phi^*\mu_2$. 
\end{theorem}

\begin{remark} {\rm It follows from (\ref{19.3})  that 
$\epsilon$ and  $\mu$ do not transform like tensors of type
(1,1). This
is due to the special role played by the underlying Euclidean
metric $g_0^{ij}=\delta^{ij}$, which is not changed
by diffeomorphisms $\Phi$. These  transformations were observed 
also in the study
of the Calder\'{o}n inverse conductivity problem.
It is shown in \cite{Sy} that, for  $\Omega\subset \R^2$,  
boundary measurements
determine the anisotropic conductivity up to same group
of transformations 
as described in Theorem \ref{group}. 
For $n\geq 3$, a similar
result is conjectured, 
based on the analysis of the linearized inverse problem, see \cite{Sy2}. 
The Calder\'{o}n problem
is closely related to the inverse problem for Maxwell's equations,
as the low-frequency limit of ${\mathcal Z}$ is related 
to the Dirichlet-to-Neumann
map for the conductivity equation \cite{La3}.}
\end{remark}

When $\epsilon$ and $\mu$ are isotropic, we obtain the
following uniqueness result.

\begin{theorem} 
\label{uniqueness in isotropic case}
Let $M\subset \R^3$ be a bounded smooth domain,
$\Gamma\subset \p M$ be open and
non-empty, $\epsilon$ and $\mu$ be smooth positive
functions on $\overline M$ and $\mathcal{Z}_{\Gamma}$
be a local admittance map for $(M,\e,\mu)$. Then $\Gamma$ and 
 ${\mathcal Z}_{\Gamma}$ determine $(M,\epsilon,\mu)$
uniquely.
\end{theorem}

Note that the knowledge of $M$ is not a priori assumed
in the above theorem.

{\bf Proof.} Assume that for $(M_1,\epsilon_1,\mu_1)$
and $(M_2,\epsilon_2,\mu_2)$ such that 
$\Gamma\subset \p M_1\cap \p M_2$ 
we have ${\mathcal Z}_{\Gamma,1}={\mathcal Z}_{\Gamma,2}$.
Then there is a diffeomorphism
 $\Phi:M_1\to M_2$ satisfying
$\Phi\big|_{\Gamma}={\rm id}$ and 
$\epsilon_1=\Phi^*\epsilon_2$, $\mu_1=\Phi^*\mu_2$. 
Since $\epsilon_1$ and $\epsilon_2$
are isotropic, it follows from the Liouville theorem  that
$\Phi$ is conformal. Since $\Phi\big|_{\Gamma}={\rm id}$,
it follows that $\Phi$ is identity.
\hfill$\Box$

\subsection{Outlook}

There are several direction to which the present work can
be extended.

1. A natural inverse  problem is the
inverse boundary spectral problem for the electric
Maxwell operator $\M_{\rm e}$. The problem is to determine the metric
$g$ and wave impedance $\alpha$, or, in the other words, $\e$ and $\mu$
from the non-zero eigenvalues $\lambda_j$ of $\M_{\rm e}$
and the normal components of the corresponding
eigenforms on $\p M$. This problem was studied  in, e.g. \cite{L2},
for the scalar
Maxwell's equations.
For
the considered anisotropic case,
  this requires significant
modifications of the method developed in this paper
and will be published elsewhere.

%

2. With the uniqueness of the inverse problem in hand, the next
issue is to study stability of the inverse problem and develop
stable reconstruction algorithms. A general approach to these questions,
in the scalar case, is introduced in \cite{K2L}, in terms of
  certain geometrical a priori bounds on $(M,g)$, with sharp results
on conditional stability in \cite{AKKLT}. Adding {\em a priori} analytical
bounds on $\alpha$, we intend to analyse these questions for anisotropic
Maxwell's equations.

\vspace{-3mm}

\subsection*{Appendix: The WKB approximation}

Denote by
$\Phi(x,y,t)=\Phi_{\lambda}(x,y,t)$ the Green's 2--form, i.e.,
the solution of
\beq
\label{11.01}
& &(\p_t^2   + \Delta_\alpha^2)\Phi_{\lambda}(x) =0
\quad \hbox{in} \,\, M\times \R_+,\\
& & \nonumber
\Phi_{\lambda} (x)|_{t=0} = 0,
\quad \p_t \Phi_{\lambda}(x)|_{t=0} = \lambda  \ud _y (x),
\quad  \Phi_{\lambda} (x)|_{\p M \times \R_+} = 0,
\eeq
where $\lambda \in \Lambda^2T^*_y M$.
Let
$B_y(\rho), \, \rho < \tau(y,\pM)$ be a domain of normal coordinates
centered at $y$, so that
\beq
\label{11.05}
g^{ij}(0)= \delta^{ij}, \, \p_kg^{ij}(0)=0.
\eeq
Rewriting equations (\ref{11.01}), componentwise,  in
these coordinates and using the unit propagation speed,
we can, instead of
(\ref{11.01}), consider the fundamental solution, $\Phi(x,y,t)$, $t< \rho$,
\beq
\label{11.06}
& &\left\{ (\p_t^2 - g^{ij} \p_i \p_j)I + B^i \p_i +C \right \}
\Phi = 0,
\quad \hbox{in $M\times ]0,\rho[$},\\
& & \nonumber
\Phi|_{t=0} = 0, \quad\p_t\Phi|_{t=0} =  I \ud (x),
\eeq
where $I$ is the $3 \times 3$ identity matrix and
   $B^i(x), \, C(x)$ are smooth
$3 \times 3$ matrices.

Following \cite{Co,Ba},
which deal with the scalar case,
we search for the solution to (\ref{11.06}) in the WKB form:
\beq
\label{11.03}
\Phi(x,t) \approx G_0(x) \, \ud (t^2 - \tau^2) +
\sum_{\ell \geq 1} G_\ell(x) \, (t^2 - \tau^2)_+^{\ell-1}/(\ell-1)!,
\eeq
where $\tau(x,y) =|x|$.
Substitution of  (\ref{11.03}) into equation (\ref{11.06})
gives rise to the recurrent system of transport equations.
The principal one
is the equation for $G_0$,
\ba
4 \tau \frac{d G_0}{d \tau}(\tau \hat{x}) +
\left \{(g^{ij} (\tau \hat{x}) \, \p_i\p_j \tau ^2 - 6)\,I
   + B^i(\tau \hat{x}) \, \p_i \tau^2 \right \}
   G_0(\tau \hat{x}) = 0,
\ea
where $\hat{x} = x/\tau$. To satisfy initial
conditions in (\ref{11.06}),
   we require that
$
G_0(0) = (2 \pi)^{-1} I.
$
By (\ref{11.05}), $g^{ij} \p_i\p_j \tau ^2\big|_{x=0} - 6=0 $. Also,
$\p_i \tau^2 \big|_{x=0} =0$. Therefore,
\bfo
\frac{1}{4 \tau}\left \{(g^{ij} (\tau \hat{x}) \, \p_i\p_j \tau
^2 - 6)\,I
   + B^i(\tau \hat{x}) \, \p_i \tau^2  \right \}
\efo
is a  smooth function of $(\tau, \hat{x})$, so that
   $G_0(x)$ is a smooth $3 \times 3$ matrix-function
  of $(\tau,\hat{x})$, for $\tau >0$.
Moreover, it can be shown that
$G_0(x)$ is also smooth at  $x=0$.

For
$G_\ell$, $\ell \geq 1$, we obtain  transport equations
\ba
& &4 \tau \frac{d G_\ell}{d \tau} +
\left \{(4\ell -6 + g^{ij}(x) \p_i\p_j \tau ^2(x)) \,I
+ B^i(x) \p_i \tau^2 \right \} G_\ell\\
& & =
   \left[   g^{ij} \p_i \p_jI - B^i \p_i -C \right ]
G_{\ell-1},
\ea
with $G_\ell(0)=0$.
If we write $G_\ell = G_0 F_\ell$, we obtain for  $F_\ell$ the equations
\beq
\label{11.07}
4 \tau \frac{d F_\ell}{d \tau} +4\ell\,F_\ell =
G_0^{-1} \, \left[   g^{ij} \p_i \p_jI - B^i \p_i -C \right ]
G_{\ell-1}, \quad F_\ell(0) =0.
\eeq
Solving the above equation, we find
\beq
\label{11.08}
F_\ell(x) = \frac14 \tau ^{-\ell} \int _0^{\tau}
G_0^{-1}(s \hat{x}) \, \left\{
\left[ g^{ij} \p_i \p_jI - B^i \p_i -C \right ]
G_{\ell-1} \right \}(s \hat{x}) \, s^{\ell-1} d s,
\eeq
which are smooth functions of $x$.
As (\ref{11.06}) is a hyperbolic system,  the right side
of (\ref{11.03}) is the asymptotics,
with respect to smoothness, of $\Phi(x,y,t)$,
when $t < \rho$.

Clearly, the asymptotic expansion (\ref{11.03}) implies decomposition
(\ref{7a1}), (\ref{remainder}).

{\bf Acknowledgements:} We would like to give our warmest thanks to
professor Alexander Katchalov for numerous useful discussions.
His lectures on non-stationary Gaussian beams
\cite{Ka1,Ka2} at Helsinki
University of Technology were paramount for our understanding
of the subject.
Also, we thank Helsinki University of Technology, Loughborough
University and University of Helsinki for their kind hospitality
and financial support.
Furthermore, the financial support of the Academy of Finland,
project 102175,
Royal society, and
EPSRC, GR/R935821/01
 is acknowledged.


\vspace{-.5cm}



\begin{thebibliography}{99}

{

\bibitem{AKKLT}
Anderson, M., Katsuda, A., Kurylev, Y., Lassas, M., Taylor M. 
Geometric Convergence, and Gel'fand's Inverse Boundary Problem.
Inventiones Mathem., {\bf 158} (2004),  261-321. 



\bibitem{Birman}
Birman, M., Solomyak, M.
The selfadjoint Maxwell operator in arbitrary domains. (Russian)
{\it Algebra i Analiz} {\bf 1} (1989), no. 1, 96--110.

\bibitem {Ba} Babich  V.M. The Hadamard ansatz,
its analogues, generalisations and applications.
{\it St. Petersburg Math. J.}, Vol. 3, No. {\bf 5} (1992), 937-972.

\bibitem {BeBl}
Belishev, M., Blagove\v s\v censkii, A.S. Direct method to solve
a nonstationary inverse problem for the wave equation.
{\it Ill-Posed Problems of Math. Phys and Anal.} (Russian)
  1988, 43-49.



\bibitem{Be1}
    Belishev, M. An approach to multidimensional inverse problems
for the wave equation.
(Russian) {\it Dokl. Akad. Nauk SSSR} {\bf 297} (1987), 524--527.

\bibitem{BeIsPSh}
    Belishev, M., Isakov, V., Pestov, L., Sharafutdinov, V.
On the reconstruction of a metric
from external electromagnetic measurements. (Russian)
{\it Dokl. Akad.
Nauk} {\bf 372}
(2000), no. 3, 298--300.

\bibitem{BeKu3}
    Belishev, M., Kurylev, Y. To the reconstruction of a Riemannian manifold
via its
spectral data (BC-method). {\it  Comm. Part. Diff. Equations} {\bf 17} (1992),
  767--804.


\bibitem{Be3}
Belishev, M., Glasman, A. A dynamic inverse problem for the Maxwell 
system: reconstruction of the velocity in the regular zone (the 
BC-method).
{\it  Algebra i Analiz} {\bf  12} (2000), 131--187; transl. St.
   Petersburg Math. J. {\bf 12} (2001),  279--316.

\bibitem{Be-Gl2}  Belishev M., Glasman A. Boundary control of the
Maxwell dynamical system: lack of controllability by topological reasons.
Mathematical and numerical aspects of wave propagation WAVES 2003, 177-182,
Springer, Berlin, 2003.

\bibitem{BI}
Belishev, M., Isakov, V. On the uniqueness of the reconstruction of 
the parameters of the Maxwell system from dynamic boundary data. 
(Russian) Zap. Nauchn. Sem. POMI {\bf 285}
   (2002), 15--32.



\bibitem{bossavit1} Bossavit,  A.
{\em \'{E}lectromagn\`{e}tisme, en vue de la mod\`{e}lisation.}
Math\'{e}matiques \& Applications 14, Springer-Verlag,
1993, xiv+174 pp.

\bibitem{bossavit2} Bossavit,  A. {\em Computational electromagnetism.
Variational formulations, complementarity, edge elements.}
  Academic Press Inc., San Diego, CA, 1998.


\bibitem {Cl}
Calder\'{o}n, A.-P. On an inverse boundary
value problem. {\it Seminar on Numerical
Analysis and
its Applications to Continuum Physics (Rio de Janeiro, 1980)}, pp. 65--73,
Soc. Brasil. Mat., Rio de Janeiro, 1980.

\bibitem{CP}
Colton,\,D., P\"aiv\"arinta,\,L.
The uniqueness of a solution to an inverse scattering problem for
electromagnetic waves.\,{\it Arch. Rational Mech. Anal.}
{\bf 119}\,(1992),\,59-70.

\bibitem{Co} Courant,  R., Hilbert, D. {\it
Methods of Mathematical Physics.
Vol.II: Partial Differential Equations},  Interscience, 1962.


\bibitem{EIsNkTa}
Eller, M., Isakov, V. Nakamura, G., Tataru, D. Uniqueness and
stability in the Cauchy problem for Maxwell's and elasticity
systems. {\it Nonlinear PDE and their applications. College de France
   Seminar}, {\bf XIV}(2002), 329--349, 2002.

\bibitem{Frankel}
Frankel, T.
{\em The geometry of physics.}
Cambridge University Press, 1997. 654 pp.

\bibitem{Ge}
Gel' fand, I.M. Some aspects of functional analysis and algebra.
{\it Proc. Intern. Cong. Math.}, {\bf 1} (1954) (Amsterdam 1957), 253--277.

\bibitem{GU}
Greenleaf, A., Uhlmann, G.
Local uniqueness for the Dirichlet-to-Neumann map via the two-plane transform.
{\it Duke Math. J.} {\bf 108} (2001), 599--617.


\bibitem{Isak1} Isakov, V. Carleman type estimates and their applications.
In: K. Binghham, Y. Kurylev, E. Somersalo (eds.): {\em New Analytical
and Geometric Methods in Inverse Problems}, Springer Verlag,
93--126.

\bibitem{IU}
Isozaki, H., Uhlmann, G.:
Hyperbolic Geometry and the Dirichlet-to-Neumann map. 
To appear in {\it Advances in Math.} 

\bibitem{Ka1}
Kachalov, A.P. Gaussian beams for the Maxwell equations on a 
manifold. (Russian) {\it Zap. Nauchn. Sem. POMI} {\bf 285} (2002), 
58--87.

\bibitem{Ka2}  Kachalov, A.P. Nonstationary
electromagnetic Gaussian beams in a nonhomogeneous anisotropic medium.
{\it Zap. Nauchn. Sem. POMI} 264 (2000), 83--100, 323.



\bibitem{KK2}
Kachalov, A., Kurylev, Y.
Multidimensional inverse problem with incomplete boundary
spectral data.
{\it Comm. Part.  Diff. Equations}, {\bf 23} (1998), 55-95.

\bibitem{KKL} Katchalov, A., Kurylev, Y., Lassas, M. {\it Inverse Boundary
Spectral Problems},  Chapman Hall/CRC, Pure and Applied Mathematics, 123.
(2001), 290 pp.


\bibitem{KKL2}
  Katchalov A., Kurylev Y., Lassas, M.
Energy measurements and equivalence of boundary
data for inverse problems on non-compact manifolds.
In {\it Geometric Methods in Inverse Problems and PDE Control} Ed. C. Croke,
I. Lasiecka, G. Uhlmann, M. Vogelius, 
 IMA volumes in Mathematics and Applications  Volume 137 (2003)
pp. 183--214. 


\bibitem{K2L}
Katsuda A., Kurylev Y., Lassas, M.
Stability and reconstruction in Gel'fand inverse boundary spectral problem.
In: K. Binghham, Y.  Kurylev, E. Somersalo (eds.): {\em New Analytical
and Geometric Methods in Inverse Problems}, 309--320, Springer-Verlag, 2003.

\bibitem{KSU}
Kenig C., Sjoestrand  J., Uhlmann  G.: 
The Calder\'on problem with partial data, preprint 
arXiv math.AP/0405486.


\bibitem{Ku1}  Kurylev, Y. Admissible groups of
transformations that preserve the boundary spectral data
in multidimensional inverse problems. (Russian) {\it Dokl. Akad. Nauk},
{\bf 327} (1992), 322--325; translated in {\it Soviet Phys. Dokl.}
{\bf 37} (1993), 544--545.

\bibitem{Ku2}  Kurylev, Y.  A multidimensional Gel'fand--Levitan inverse
boundary problem {\it Diff. Equat. and Math. Phys.} (Birmingham, AL, 1994),
117-131, Int. Press, 1995.


\bibitem{Ku3} Kurylev, Y.
An inverse boundary problem for the Schr\"odinger operator with
magnetic field.
{\it J. Math. Phys.} {\bf 36} (1995), no. 6, 2761--2776.

\bibitem{Ku5}
    Kurylev, Y.  Multidimensional Gel'fand inverse problem
and boundary distance map, {\it Inverse Problems Related with
Geometry} (ed. H.Soga) (1997), 1-15.

\bibitem{KL1}
Kurylev, Y., Lassas, M.
    Gelf'and inverse problem for a quadratic operator pencil.
{\it J.
Funct. Anal.} {\bf 176} (2000), no. 2, 247--263.

\bibitem{KL2}
Kurylev, Y., Lassas, M. Hyperbolic inverse problem with data on a part of
the boundary.
{\it Differential Equations and Mathematical Physics (Birmingham, AL, 1999)},
259--272, AMS/IP Stud. Adv. Math., {\bf 16},
Amer. Math. Soc., 2000.


\bibitem{KL3}  Kurylev, Y., Lassas, M. Hyperbolic inverse
problem and unique continuation of Cauchy data of solutions along the
boundary, {\it Proc. Roy. Soc. Edinburgh, Ser. A},
  {\bf 132} (2002), no. 4, 931--949.


\bibitem{KLS}  Kurylev, Y., Lassas, M., Somersalo  E.
  Reconstruction of a manifold-from electromagnetic boundary 
measurements, {\it Contemporary Mathematics}.
333(2002), Inverse Problems: Theory and Applications. 
Ed. G. Alessandrini, G. Uhlmann, 147--162. 

\bibitem{KLS2} 
 Kurylev, Y., Lassas, M., Somersalo  E. 
Focusing waves in elctromagnetic inverse problems,
{\it Contemporary Mathematics} 348(2004), 11-22, Inverse Problems and Spectral Theory Ed. H. Isozaki.





\bibitem{LanL}
Landau, L., Lifshitz, E. {\it Course of theoretical physics, Vol. 2.
The classical theory of fields.} Pergamon Press, 1975. xiv+402 pp.



\bibitem{langer} Langer,  R.E. An inverse problem in differential
equations. {\it Bull. Amer. Math. Soc} {\bf 39} (1933), 814--820.

\bibitem{LTr} Lasiecka, I., Triggiani, R.
{\it Control Theory for Partial Differential Equations.}
  Cambridge Univ. Press, 2000, pp. 645--1067.


\bibitem{L2}Lassas,
M.  Inverse boundary spectral problem for non-selfadjoint Maxwell's
equations with incomplete data.
{\it Comm. Part. Diff. Equations}
{\bf 23} (1998), 629--648.

\bibitem{La3}
  Lassas, M. The impedance imaging problem as a low-frequency limit.
{\it Inverse Problems} {\bf 13} (1997), no. 6, 1503--1518.



\bibitem{LTU}
Lassas, M. Taylor, M., Uhlmann, G.
The Dirichlet-to-Neumann map for complete Riemannian manifolds with
boundary,
{\it Comm. Anal. Geom.} {\bf 11} (2003), 207-22.

\bibitem{LU}
Lassas, M., Uhlmann, G. Determining Riemannian manifold
from boundary measurements.
{\it Ann. Sci. Ecole Norm.
Sup.} (4) {\bf 34} (2001), no. 5, 771--787.

\bibitem{LeU}
Lee, J., Uhlmann, G.
Determining anisotropic real-analytic conductivities by boundary
measurements.
{\it Comm. Pure Appl. Math.} {\bf 42} (1989), no. 8, 1097--1112.






\bibitem{Na1}
Nachman, A.
    Reconstructions from boundary measurements.
{\it  Ann. of Math.} (2) {\bf
128} (1988),
no. 3, 531--576.

\bibitem{Na2}
Nachman, A. Global uniqueness for a two-dimensional
    inverse boundary value problem. {\it Ann.
of Math.} (2) {\bf 143} (1996), no. 1, 71--96.

\bibitem{NaSyU}
Nachman, A., Sylvester, J. Uhlmann, G. An $n$-dimensional Borg-Levinson
theorem. {\it Comm. Math. Phys.} {\bf 115 }(1988), no. 4, 595--605.

\bibitem{NvKh}
Novikov, R., Khenkin, G.
The $\overline\partial$-equation in the multidimensional inverse
scattering problem. (Russian)
{\it Uspekhi Mat. Nauk} {\bf 42} (1987), no. 3, 93--152.


\bibitem {OPS}
Ola, P., P\"aiv\"arinta, L., Somersalo, E.
    An inverse boundary value problem in
electrodynamics. {\it Duke Math. J.} {\bf 70} (1993), no. 3, 617--653.

\bibitem{OS}
    Ola, P., Somersalo, E.
Electromagnetic inverse problems and generalized Sommerfeld
potentials. {\it SIAM J. Appl. Math.} {\bf 56} (1996), no. 4, 1129--1145.

\bibitem{paquet}
Paquet, L.
Mixed problems for the Maxwell system  (French)
{\it Ann. Fac. Sci. Toulouse Math.} (5) {\bf 4} (1982), no. 2, 103--141.

\bibitem{picard} Picard,  R. On the low frequency asymptotics
in electromagnetic theory. {\it J. Reine Angew. Math.} {\bf 394} 
(1984), 50--73.

\bibitem{rom}
Romanov, V.G. An inverse problem of electrodynamics, {\it Dokl.
Mathem., }  {\bf 66} (2002), no. 2, 200 -- 205.

\bibitem{Ru} Russell, D. Controllability and stabilizability theory 
for linear partial differential equations: recent progress and open 
questions.
{\it  SIAM Rev.} {\bf 20} (1978), 639--739.

\bibitem{schlichter} Schlichter,  L.B. An inverse boundary value problem
in electrodynamics. {\it Physics} {\bf 4} (1953), 411--418.

\bibitem{Sc}  Schwarz, G.
{\it Hodge decomposition--a method for solving boundary value problems}.
  Springer-Verlag,  1995. 155 pp.


\bibitem{SIC} Somersalo,  E.,  Isaacson D., Cheney  M. A linearized inverse
boundary value problem for Maxwell's equations. {\it J. Comp. Appl. Math.}
{\bf 42} (1992), 123--136.

\bibitem{Sy}
Sylvester, J. An anisotropic inverse boundary value problem.
{\it Comm. Pure Appl. Math.} 43 (1990), no. 2, 201--232

\bibitem{Sy2} Sylvester, J. Linearizations of anisotropic inverse problems.
In: L. P\"aiv\"arinta and E. Somersalo (eds.): {\em Inverse Problems
in Mathematical Physics}. Lecture Notes in Physics 422, Springer-Verlag, 1993.

\bibitem{SyU}
    Sylvester, J., Uhlmann, G. A global uniqueness theorem for an inverse
boundary value
problem. {\it Ann. of Math.} (2) {\bf 125} (1987), no. 1, 153--169.


\bibitem{Ta1}
    Tataru, D.
    Unique continuation for solutions to PDEs; between H\"ormander's theorem
and
Holmgren's theorem. {\it Comm. Part. Diff. Equations}
{\bf 20} (1995), no.
5-6, 855--884.

\bibitem{Ta3}
Tataru, D. Unique continuation for operators
    with partially analytic coefficients. {\it J. Math.
Pures Appl.} (9) {\bf 78} (1999), no. 5, 505--521.

\bibitem{whitney} Whitney,  H. {\em Geometric integration theory.}
Princeton University Press, 1957.

}
\end{thebibliography}
\end{document}